\documentclass[12pt, reqno]{amsart}
\usepackage{graphicx}
\usepackage{tikz}

\usepackage{subcaption}

\usetikzlibrary{shapes, positioning,patterns}
\usetikzlibrary{trees}
\usetikzlibrary{fit}
\usepackage{float}
\usepackage{ulem}
\usepackage[justification=centering]{caption}
\usepackage[all]{xy}
\usepackage{amsmath}
\usepackage{amsfonts}
\usepackage{amssymb}
\usepackage[pagebackref]{hyperref}
\hypersetup{backref, pagebackref, colorlinks=true}
\usepackage{cite}
\usepackage{color}

\usepackage{epic}

\usepackage{lmodern}

\newtheorem{thm}{Theorem}[section]
\newtheorem{defi}{Definition}[section]

\newtheorem{prop}{Proposition}[section]
\newtheorem{cor}{Corollary}[section]

\newtheorem{remark}{Remark}

\usepackage{amsthm,amsmath,amssymb,url,cite,color}

\numberwithin{equation}{section}

 \usepackage{times}

\begin{document}

\title[ the Narayana polynomials and the Motzkin polynomials]{New 
 Refinements of Narayana polynomials \\[5pt]
 and Motzkin polynomials}

\author[J.J.W. Dong]{Janet J.W. Dong}
\address[Janet J.W. Dong]{Department of Mathematics, Shaoxing University, Shaoxing 312000, P.R. China}
\email{dongjinwei@tju.edu.cn}

\author[L.R. Du]{Lora R. Du}
\address[Lora R. Du]{Center for Applied Mathematics and KL-AAGDM, Tianjin University, Tianjin 300072, P.R. China}
\email{loradu@tju.edu.cn}

\author[K.Q. Ji]{Kathy Q. Ji}
\address[Kathy Q. Ji]{Center for Applied Mathematics and KL-AAGDM, Tianjin University, Tianjin 300072, P.R. China}
\email{kathyji@tju.edu.cn}

\author[D.T.X. Zhang]{Dax T.X. Zhang}
\address[Dax T.X. Zhang]{College of Mathematical Science \& Institute of Mathematics and Interdisciplinary Sciences, Tianjin Normal University, Tianjin 300387, P. R. China}
\email{zhangtianxing6@tju.edu.cn}

 \date{\today}
\begin{abstract} Chen, Deutsch and Elizalde  introduced a refinement of the Narayana polynomials by distinguishing between old (leftmost child) and young leaves of plane trees. They also provided a refinement of Coker's formula by constructing a bijection. In fact,  Coker's formula establishes a connection between the Narayana polynomials and the Motzkin polynomials, which implies the $\gamma$-positivity of the Narayana polynomials.  In this paper, we introduce the polynomial $G_{n}(x_{11},x_{12},x_2;y_{11},y_{12},y_2)$, which further refine the Narayana polynomials by considering  leaves of plane trees that have no siblings. We obtain the generating function for $G_{n}(x_{11},x_{12}, x_2;y_{11},y_{12},y_2)$.  To achieve further refinement of Coker's formula based on the polynomial $G_n(x_{11},x_{12},x_2;y_{11},y_{12},y_2)$, we consider a refinement $M_n(u_1,u_2,u_3;v_1,v_2)$ of the Motzkin polynomials by classifying the old leaves of a tip-augmented plane tree into three categories and the young leaves into two categories. The generating function for $M_n(u_1,u_2,u_3;v_1,v_2)$ is  also established,  and the refinement of Coker's formula is immediately derived by combining the generating function   for $G_n(x_{11},x_{12},x_2;y_{11},y_{12},y_2)$ and the generating function   for $M_n(u_1,u_2,u_3;v_1,v_2)$. 
We derive several interesting consequences from this  refinement of Coker's formula,   including the new symmetries of vertices of  plane trees, Euler transformation of the Narayana polynomials and the Motzkin polynomials due to Lin-Kim  and  the real-rootedness of the Motzkin polynomials. The method used in this paper is the grammatical approach introduced by Chen. We develop a unified  grammatical approach to exploring polynomials associated with the statistics defined on  plane trees. As you will see, the derivations of  the generating functions for $G_n(x_{11},x_{12},x_2;{y}_{11},{{y}}_{12},y_2)$ and  $M_n(u_1,u_2,u_3;v_1,v_2)$ become quite simple once their grammars are established.
 
\end{abstract}

\maketitle

\section{Introduction}
The main objective of this paper is to present a grammatical approach to studying the Catalan numbers, the Narayana numbers,  the Motzkin numbers and related topics.
The Catalan numbers and the Motzkin numbers are  classical topics in enumerative combinatorics, which have been extensively investigated
for decades, see \cite{Donaghey-1977,  Donaghey-Shapiro-1977, Petersen-2015, Stanley-2015}.  

The Narayana numbers can be viewed as a refinement of  the Catalan numbers. It is well-known that  the Catalan numbers $C_n=\frac{1}{n+1}\binom{2n}{n}$  counts the number of unlabelled (rooted) plane trees with $n$ edges. In this context,    the Narayana numbers $N(n,k)$ provide a more detailed enumeration 
by counting  the number of such plane trees with $n$ edges and $k$ leaves. 
The explicit formula for $N(n,k)$ is given by 
\begin{equation}\label{exp-Na}
N(n,k)=\frac{1}{n}\binom{n}{k}\binom{n}{k-1} 
\end{equation}
with $N(0,0)=1$ and $N(0,k)=0$ for $k\geq 1$. The Narayana polynomials are defined as  
\[N_n(x)=\sum_{k=1}^{n}N(n,k)x^k\]
and the generating function for the Narayana polynomial $N_n(x)$ is given by 
\begin{equation} \label{gf-Narayana}    \sum_{n=0}^{\infty}N_n(x)q^n=\frac{1+q-xq-\sqrt{1-2(1+x)q+(x-1)^2q^2}}{2q},
\end{equation}
see Deutsch \cite{Deutsch-1999} and  Petersen \cite[p. 24]{Petersen-2015}.

Let $\mathcal{P}_n$ denote the set of unlabelled (rooted) plane trees with $n$ edges. For a plane tree $T \in \mathcal{P}_n$, let ${\rm leaf}(T)$ denote the number of leaves in $T$. The Narayana polynomials can be expressed as 
\begin{equation}\label{com-narayana}
N_{n}(x)=\sum_{T \in \mathcal{P}_n}x^{{\rm leaf}(T)}.
\end{equation}

\begin{figure}[H]
\centering
\begin{tikzpicture}[scale=0.5]
	    \draw[fill](0,0)circle(2pt)--(0,-1)circle(2pt)--(0,-2)circle(2pt)--(0,-3)circle(2pt)node[below]{$x$};
\end{tikzpicture}\qquad
\begin{tikzpicture}[scale=0.5]
\draw[fill](0,0)circle(2pt)--(0,-1)circle(2pt)--(-1,-2)circle(2pt)node[below]{$x$};
\draw[fill](0,-1)--(1,-2)circle(2pt)node[below]{$x$};
\end{tikzpicture}\qquad
\begin{tikzpicture}[scale=0.5]
\draw[fill](0,0)circle(2pt)--(-1,-1)circle(2pt)--(-1,-2)circle(2pt)node[below]{$x$};
\draw[fill](0,0)--(1,-1)circle(2pt)node[below]{$x$};
\end{tikzpicture}\qquad
\begin{tikzpicture}[scale=0.5]
\draw[fill](0,0)circle(2pt)--(1,-1)circle(2pt)--(1,-2)circle(2pt)node[below]{$x$};
\draw[fill](0,0)--(-1,-1)circle(2pt)node[below]{$x$};
\end{tikzpicture}\qquad
\begin{tikzpicture}[scale=0.5]
\draw[fill](0,0)circle(2pt)--(0,-1)circle(2pt)node[below]{$x$};
\draw[fill](0,0)--(1,-1)circle(2pt)node[below]{$x$};
\draw[fill](0,0)--(-1,-1)circle(2pt)node[below]{$x$};
\end{tikzpicture}
\caption{$3$-edge plane trees.}
    \label{example-plane trees-1}
\end{figure}
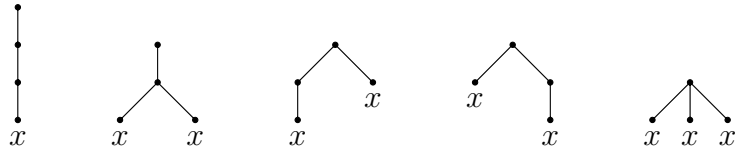

For example, as illustrated by the plane trees in Fig. \ref{example-plane trees-1}, we observe that $N_3(x)=x+3x^2+x^3$. 

The Narayana polynomials have been studied by several authors, see, for example,  Bonin, Shapiro and Simion \cite{Bonin-Shapiro-Simion-1993}, Coker \cite{Coker-2003}, 
Mansour and Sun \cite{Mansour-Sun-2008, Mansour-Sun-2009}, Rogers \cite{Rogers-1981}, Rogers  and  Shapiro \cite{Rogers-Shapiro-1981}  and Sulank \cite{Sulanke-2000, Sulanke-2004}.  
In particular,  the Narayana polynomial $N_n(x)$ is also the $h$-polynomial of the simplicial complex dual to an associahedron of type $A_n$, see Formin and Reading \cite{Fomin-Reading-2007}.  

Using \eqref{exp-Na}, we see that $N(n,k)=N(n, n-k+1)$ for $1\leq k\leq n$, leading to  the following combinatorial consequence.  Elegant combinatorial proofs of this proposition were provided by Chen \cite{Chen-1990, Chen-1999}  and Schmitt and Waterman \cite{Schmitt-Waterman-1994}.
\begin{prop} \label{sym-Narayana}
The number of  plane trees with $n$ edges and $k$ leaves is equal to the number of plane trees with $n$ edges and $k$ interior vertices. 
\end{prop}
It follows that $N_n(x)$ is a symmetric polynomial. Recall that   a polynomial $a_0+a_1q+\cdots+a_Nq^n$
  with integer coefficients is called symmetric if $a_j=a_{N-j}$ for $0\leq j\leq N$. 
 It is known that the Narayana polynomial $N_n(x)$ has only distinct non-positive real zeros, see Petersen \cite[Problem 4.7]{Petersen-2015}. As pointed out by Br\"{a}nd\'{e}n\cite[Remark 7.3.1]{Branden-2015}, the real-rootedness of a polynomial with symmetric and non-negative coefficients implies the $\gamma$-positivity of
such polynomial. In fact, the $\gamma$-positivity of the Narayana polynomial $N_n(x)$ can be deduced from the following formula due to Coker \cite{Coker-2003}. 

\begin{thm}[Coker]\label{cokerthm} For $n\geq 1$, 
\begin{equation} \label{eq-coker}
N_{n}(x)=\sum_{k=1}^{\lfloor \frac{n+1}{2} \rfloor } \binom{n-1}{2k-2}C_{k-1} x^k (1+x)^{n-2k+1}. 
\end{equation} 
\end{thm}
 Coker \cite{Coker-2003} obtained this formula by using the Lagrange inversion formula. Chen, Yan and Yang \cite{Chen-Yan-Yang-2008} provided a bijective proof of this formula by building a  bijection between Dyck paths and $2$-Motzkin paths. 
By comparing the coefficients of $x^n$ in the formula \eqref{eq-coker}, one can recover an identity of Simion and Ulman \cite{Simion-Ullman-1991}  expressing the Narayana numbers by the Catalan numbers.  The bijective proof of the identity of Simion and Ulman  was given by Chen, Deng and Du \cite{Chen-Deng-Du-2005}.

Let ${P}_e(n)$ (${P}_o(n)$) denote the number of plane trees with $n$ edges and an even (odd) number of leaves.  Setting $x=-1$ in \eqref{eq-coker} and using \eqref{com-narayana},  we obtain 
\begin{cor} For $n\geq 1$, 
\begin{align*}
{P}_e(2n)-{P}_o(2n) &=0,\\[5pt]
{P}_e(2n+1)-{P}_o(2n+1) &=(-1)^{n+1} C_n. 
\end{align*}
\end{cor}
These relations were also obtained by Bonin, Shapiro and Simion \cite{Bonin-Shapiro-Simion-1993}, Eu, Liu and Yeh \cite{Eu-Liu-Yeh-2004}  and Klazar \cite{Klazar-2003} 
 independently.  A bijective proof of these two relations was given by Chen, Shapiro and Yang \cite{Chen-Shapiro-Yang-2006}. 

In fact,   Coker's formula \eqref{eq-coker}  built a connection between the Narayana numbers and the Motzkin numbers, see Chen and Pan \cite{Chen-Pan-2017}. The  Motzkin numbers,  introduced in \cite{Motzkin-1948} and  denoted $M_n$,  are defined by the following recurrences:  
\begin{align*}
M_0=1, \quad & M_{n+1}=M_n+\sum_{k=0}^{n-1} M_k M_{n-1-k} \quad \text{for} \quad n\geq 0. 
\end{align*}
The following two identities describe the relationships between  the Motzkin numbers $M_n$ and the Catalan numbers 
$C_n$:    
\begin{equation}\label{rmotcat}
(a).\ M_n=\sum_{k=0}^{\lfloor\frac{n}{2}\rfloor} \binom{n}{2k} C_k \quad \text{and} \quad (b).\ C_{n+1}=\sum_{k=0}^{n}\binom{n}{k}M_k.
\end{equation}
Combinatorial proofs of the above two relations were given by Donaghey\cite{Donaghey-1977}. 

Chen and Pan \cite{Chen-Pan-2017} defined the following  two-variable  Motzkin polynomials:   \begin{equation}\label{Motzkinpoly}
    M_n(u;v) = \sum_{k=0}^{\lfloor \frac{n}{2}\rfloor}\binom{n}{2k}C_k u^{k+1}v^{n-2k }.
\end{equation}
It should be noted that  when $u$ and $v$ are positive integers,  the Motzkin polynomial  $M_n(u;v)$ was termed the generalized Motzkin numbers by Sun \cite{Sun-2014}. Sun \cite{Sun-2014} also obtained the generating function for $M_n(u;v)$:

\begin{thm}[Sun] \label{Sun-gene-thm} 
 \begin{equation}\label{Sun-gene}
     \sum_{n\geq 0}M_n(u;v)q^n=
 \frac{1-vq-\sqrt{1-2vq+(v^2-4u)q^2 }}{2q^2}. 
 \end{equation}
\end{thm}

In the notation of $N_n(x)$ and $M_n(u;v)$, we find that Coker's result \eqref{eq-coker} can be recast as follows: 
\begin{thm}[Coker]\label{cokerthmaa}
For $n\geq 0$, 
\begin{equation} \label{two-variable}
N_{n+1}(x)=M_{n}(x;1+x).
\end{equation}
\end{thm}
In fact, Theorem \ref{cokerthmaa} follows immediately from   comparing the generating function \eqref{gf-Narayana} for $N_n(x)$ with the generating function \eqref{Sun-gene} for $M_n(u;v)$.

Chen, Deutsch and Elizalde \cite{Chen-Deutsch-Elizalde-2006} refined the Narayana polynomials by classifying the leaves of a plane tree into old and young leaves. Specifically, a leaf is considered an old leaf if it is the leftmost child of its parent, and a young leaf otherwise. It should be noted that if a leaf is the only child of its parent, Chen, Deutsch and Elizalde also classified this leaf as an old leaf.

  Let ${\rm oleaf}(T)$ and ${\rm yleaf}(T)$  denote the numbers of old leaves and young leaves in a plane tree $T$, respectively. Chen, Deutsch and Elizalde \cite{Chen-Deutsch-Elizalde-2006} defined the following polynomials: 
  \begin{equation}\label{defi-C}
   G_n(x_1,x_2)=\sum_{T \in \mathcal{P}_n}x_1^{{\rm oleaf}(T)}x_2^{{\rm yleaf}(T)}. 
\end{equation}
When $x_1=x_2=x$ in $ G_n(x_1,x_2)$, we recover the Narayana polynomial $N_n(x)$.

\begin{figure}[H]
\centering
\begin{tikzpicture}[scale=0.5]
	    \draw[fill](0,0)circle(2pt)--(0,-1)circle(2pt)--(0,-2)circle(2pt)--(0,-3)circle(2pt)node[below]{$x_1$};
\end{tikzpicture}\qquad
\begin{tikzpicture}[scale=0.5]
\draw[fill](0,0)circle(2pt)--(0,-1)circle(2pt)--(-1,-2)circle(2pt)node[below]{$x_1$};
\draw[fill](0,-1)--(1,-2)circle(2pt)node[below]{$x_2$};
\end{tikzpicture}\qquad
\begin{tikzpicture}[scale=0.5]
\draw[fill](0,0)circle(2pt)--(-1,-1)circle(2pt)--(-1,-2)circle(2pt)node[below]{$x_1$};
\draw[fill](0,0)--(1,-1)circle(2pt)node[below]{$x_2$};
\end{tikzpicture}\qquad
\begin{tikzpicture}[scale=0.5]
\draw[fill](0,0)circle(2pt)--(1,-1)circle(2pt)--(1,-2)circle(2pt)node[below]{$x_1$};
\draw[fill](0,0)--(-1,-1)circle(2pt)node[below]{$x_1$};
\end{tikzpicture}\qquad
\begin{tikzpicture}[scale=0.5]
\draw[fill](0,0)circle(2pt)--(0,-1)circle(2pt)node[below]{$x_2$};
\draw[fill](0,0)--(1.5,-1)circle(2pt)node[below]{$x_2$};
\draw[fill](0,0)--(-1.5,-1)circle(2pt)node[below]{$x_1$};
\end{tikzpicture}
\caption{$3$-edge plane trees.}
    \label{example-plane trees-2}
\end{figure}
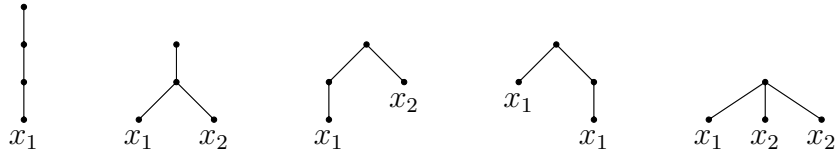

Fig. \ref{example-plane trees-2} shows five  plane trees with three edges, where the old leaves are labeled by $x_1$ and the young leaves are labeled by $x_2$. Hence, we have  
\[G_3(x_1,x_2)=x_1+2x_1x_2+x_1^2+x_1x_2^2.\]

Chen, Deutsch and Elizalde \cite{Chen-Deutsch-Elizalde-2006} derived the following generating function for $ G_n(x_1,x_2)$ by establishing an equation for the generating function using a decomposition of plane trees.  

\begin{thm}[Chen-Deutsch-Elizalde]
    \begin{equation}\label{OY-1}
  \sum_{n\geq 0} G_n(x_1,x_2)q^n= \frac{1+(1-x_2)q-\sqrt{1-2(1+x_2)q+(1-4x_1+2x_2+x_2^2)q^2}}{2q}.
    \end{equation}
\end{thm}

Let $P(n,i,j)$ denote the number of plane trees with $n$ edges, $i$ old leaves and $j$ young leaves, we see that 
\[G_n(x_1,x_2)=\sum_{i,j}P(n,i,j)x_1^i x_2^j.\]
By applying the Lagrange inversion formula to the generating function for $G_n(x_1,x_2)$, Chen, Deutsch and Elizalde obtained the following explicit formula for $P(n,i,j)$:
\begin{equation}
P(n,i,j)=\frac{1}{n}\binom{n}{i}\binom{n-i}{j}\binom{n-i-j}{i-1}.
\end{equation}
Moreover, they  established a bijection between the set of plane trees with $n$ edges and the set of 2-Motzkin paths of length $n-1$, leading to the following refinement of Coker's formula: For $n\geq 0$,  
    \begin{align} \label{Chen-Coker}
        &\sum_{i=1}^{\lfloor \frac{n+2}{2}\rfloor}\sum_{j=0}^{n+2-2i}\frac{1}{n+1}\binom{n+1}{i}\binom{n+1-i}{j}\binom{n+1-i-j}{i-1}x_1^ix_2^j \notag\\[5pt]
        &=\sum_{k=0}^{\lfloor \frac{n}{2}\rfloor} C_k\binom{n}{2k} x_1^{k+1}(1+x_2)^{n-2k}.
    \end{align}
Using the notation  $G_n(x_1, x_2)$ and $M_n(u;v)$, we can state  the identity  \eqref{Chen-Coker} as follows: 
\begin{thm}[Chen-Deutsch-Elizalde]\label{chenrel}For $n\geq 0$, 
\begin{equation}\label{eqchenrel}
G_{n+1}(x_1, x_2)=M_n(x_1;1+x_2).
\end{equation}     
\end{thm}
In fact, this relation can be immediately derived by  comparing the generating function \eqref{OY-1}  for $G_n(x_1,x_2)$ with the generating function \eqref{Sun-gene}  for $M_n(u;v)$.

In this paper, we further refine the Narayana polynomials by considering leaves without any siblings. It should be emphasized that Chen, Deutsch and Elizalde \cite{Chen-Deutsch-Elizalde-2006} classified a leaf without any siblings as an old leaf. More precisely, we define a leaf without any siblings as a singleton leaf. A leaf with siblings is considered as an elder leaf if it is the leftmost child of its parent, and a young leaf otherwise. An interior vertex is called  a young interior vertex if it is not a parent of either a singleton leaf or an elder leaf.

Let ${\rm sleaf}(T)$ and ${\rm eleaf}(T)$  denote the numbers of singleton leaves and elder leaves in a plane tree $T$, respectively. Let ${\rm sint}(T)$ and ${\rm eint}(T)$  denote the numbers of parents of singleton leaves and elder leaves in $T$ respectively, and let ${\rm yint}(T)$ denote the number of young interior vertices in $T$.
By definition, it is evident that 
\[{\rm sleaf}(T)={\rm sint}(T), \quad  {\rm eleaf}(T)={\rm eint}(T)\]
and 
for $T \in  \mathcal{P}_n$, 
\[2{\rm sleaf}(T)+2{\rm eleaf}(T)+{\rm yint}(T)=n+1.\]

\begin{figure}[H]
\centering
\begin{tikzpicture}[scale=0.5]
	    \draw[fill](0,0)circle(2pt)node[right]{$y_{2}$}--(0,-1)circle(2pt)node[right]{$y_{2}$}--(0,-2)circle(2pt)node[right]{$y_{11}$}--(0,-3)circle(2pt)node[right]{$x_{11}$};
\end{tikzpicture}\qquad
\begin{tikzpicture}[scale=0.5]
\draw[fill](0,0)circle(2pt)node[right]{$y_{2}$}--(0,-1)circle(2pt)node[right]{$y_{12}$}--(-1,-2)circle(2pt)node[below]{$x_{12}$};
\draw[fill](0,-1)--(1,-2)circle(2pt)node[below]{$x_{2}$};
\end{tikzpicture}\qquad
\begin{tikzpicture}[scale=0.5]
\draw[fill](0,0)circle(2pt)node[above]{$y_{2}$}--(-1,-1)circle(2pt)node[left]{$y_{11}$}--(-1,-2)circle(2pt)node[below]{$x_{11}$};
\draw[fill](0,0)--(1,-1)circle(2pt)node[below]{$x_{2}$};
\end{tikzpicture}\qquad
\begin{tikzpicture}[scale=0.5]
\draw[fill](0,0)circle(2pt)node[above]{$y_{12}$}--(1,-1)circle(2pt)node[right]{$y_{11}$}--(1,-2)circle(2pt)node[below]{$x_{11}$};
\draw[fill](0,0)--(-1,-1)circle(2pt)node[below]{$x_{12}$};
\end{tikzpicture}\qquad
\begin{tikzpicture}[scale=0.5]
\draw[fill](0,0)circle(2pt)node[above]{$y_{12}$}--(0,-1)circle(2pt)node[below]{$x_{2}$};
\draw[fill](0,0)--(1.5,-1)circle(2pt)node[below]{$x_{2}$};
\draw[fill](0,0)--(-1.5,-1)circle(2pt)node[below]{$x_{12}$};
\end{tikzpicture}
\caption{$3$-edge plane trees.}
    \label{example-plane trees}
\end{figure}
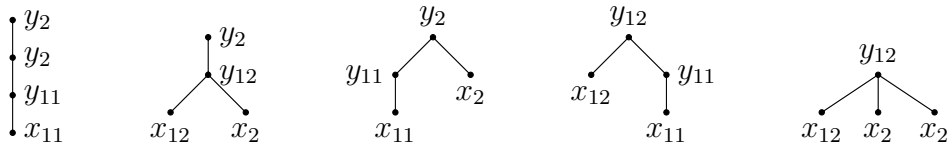

We define the following polynomials: 
\begin{defi}[New refinement of the Narayana polynomials] For $n\geq 2$, 
\begin{align}
&G_n({x}_{11},{{{x}}}_{12},x_2;{y}_{11},{{y}}_{12},y_2) \nonumber \\[5pt]
&=\sum_{T\in \mathcal{P}_{n}} {x}_{11}^{{\rm sleaf}(T)} {x}_{12}^{{\rm eleaf}(T)} x_2^{{\rm yleaf}(T)} y_{11}^{{\rm sint}(T)} {y}_{12}^{{\rm eint}(T)}{y}_2^{{\rm yint}(T)} 
\end{align}
with the convention that   
\[G_0(x_{11},x_{12},x_2;y_{11},y_{12},y_2)=y_2 \quad \text{and} \quad  G_1(x_{11},x_{12},x_2;y_{11},y_{12},y_2)=x_{12}y_{12}.\]
\end{defi}
Fig. \ref{example-plane trees} gives the list of plane trees in $\mathcal{P}_3$, where the singleton leaves are labeled by $x_{11}$, the elder leaves  are labeled by $x_{12}$, 
 the young leaves are labeled by $x_2$, the parents of singleton leaves are labeled by $y_{11}$, the parents of elder leaves are labeled by $y_{12}$, and the young interior vertices are labeled by $y_2$.  It follows that 
\begin{align*}
&G_3({x}_{11},x_{12},x_2;{y}_{11},{{y}}_{12},y_2)\\[5pt]
&=x_{11}y_{11}y_2^2+x_{12}x_2y_{12}y_2+x_{11}x_2y_{11}y_2+x_{11}x_{12}y_{11}y_{12}+x_{12}x_2^2y_{12},
\end{align*}
which is a homogeneous polynomial of degree four. 

Setting $x_{11} =x_{12} = x_1$ and $y_{11}=y_{12}=y_2=1$ in $G_n({x}_{11},{{{x}}}_{12},x_2;{y}_{11},{{y}}_{12},y_2)$ yields the polynomial $G_n(x_1, x_2)$ as introduced by Chen, Deutsch and Elizalde. Setting $x_{11} =x_{12} = x_2=x$ and $y_{11}=y_{12}=y_2=1$ in $G_n({x}_{11},{{{x}}}_{12},x_2;{y}_{11},{{y}}_{12},y_2)$ recovers the Narayana polynomial $N_n(x)$.

\begin{thm}\label{thm-gf-abxy} We have
{ \begin{align}\label{gf-abxy}
&\sum_{n\geq 0}G_{n}(x_{11},x_{12},x_2;y_{11},y_{12},y_2)q^n \nonumber \\[5pt]
&= \frac{1+a_0q+a_1q^2-\sqrt{ 1+a_2q+a_3q^2 + a_4q^3+a_5q^4}}{2q},
\end{align}}
where  
\begin{align*}
a_0&=y_2-x_2,\\[2pt]
a_1&=x_{12}y_{12}-x_{11}y_{11},\\[2pt]
a_2&=-2(x_2+y_2),\\[2pt]
a_3&=(x_2+y_2)^2-2(x_{11}y_{11}+x_{12}y_{12}),\\[2pt]
a_4&=2(x_{11}y_{11}-x_{12}y_{12})(x_2-y_2),\\[2pt]
a_5&=(x_{11}y_{11}-x_{12}y_{12})^2.
\end{align*} 

\end{thm}
From \eqref{gf-abxy}, we see that for $n\geq 2$, 
    $$G_{n}(x_{11},x_{12},x_2;y_{11},y_{12},y_2)=G_{n}(x_{12},x_{11},y_2;y_{12},y_{11},x_2).$$ 
More precisely, we have the following combinatorial assertion: 
\begin{prop}
When $n\geq 2$, the number of plane trees with $n$ edges, $k$ singleton leaves, $l$ elder leaves, $i$ young leaves and $j$ young interior vertices equals the number of plane trees with $n$ edges, $l$ singleton leaves, $k$ elder leaves, $j$ young leaves and $i$ young interior vertices.  
\end{prop}
In particular, we have the  following consequence. 
\begin{cor} \label{sym-old-young}
{\rm (a)}   When $n\geq 2$, the number of plane trees with $n$ edges, $k$ singleton leaves, $l$ elder leaves equals the number of plane trees with $n$ edges, $l$ singleton leaves, $k$ elder leaves. 

{\rm (b)} When $n\geq 2$, the number of plane trees with $n$ edges, $i$ young leaves and $j$ { young interior vertices}  equals the number of plane trees with $n$ edges, $j$ young leaves and $i$ { young interior vertices}. 
\end{cor}
It should be noted that Corollary \ref{sym-old-young}  (b) can be viewed as the refinement of the symmetry of the Narayana polynomials (Proposition \ref{sym-Narayana}).

Setting $x_{11}=x_{12}=x_{1}$ and $y_{11}=y_{12}=y_2=1$  in \eqref{gf-abxy},  we recover the generating function \eqref{OY-1} for $ G_n(x_1,x_2)$. In fact, 
we can consider the following four-variable polynomials: for $n\geq 1$, 
\begin{equation} \label{eq-fourvar}
   G_n(x_1,x_2; y_1,y_2)=\sum_{T \in \mathcal{P}_n}x_1^{{\rm oleaf}(T)}x_2^{{\rm yleaf}(T)}y_1^{{\rm oint}(T)}y_2^{{\rm yint}(T)}, 
\end{equation}
where ${\rm oint}(T)$  counts the number of interior vertices, i.e., the vertices that are the  parents of old leaves. As a preliminary condition, we set $G_0(x_{1},x_2;y_{1},y_2)=y_2.$
It is evident that 
\[{\rm oint}(T)={\rm eint} (T)+{\rm sint}(T),\]
and so for $n\geq 0$, 
  \begin{align*}
G_n(x_{1},x_{1},x_2;y_{1},y_{1},y_2)&=G_n(x_1,x_2; y_1,y_2).
 \end{align*}
Thus,  the  generating function for  $G_n(x_1,x_2; y_1,y_2)$ is derived by  setting $x_{11}=x_{12}=x_1$ and $y_{11}=y_{12}=y_1$ in \eqref{gf-abxy}. 
\begin{thm}\label{gf-Chen-D-f} We have
    \[\sum_{n\geq 0}G_n(x_1,x_2;y_1,y_2)q^n=\frac{1+(y_2-x_2)q-\sqrt{1-2(x_2+y_2)q+((x_2+y_2)^2-4x_1y_1)q^2}}{2q}.\]
\end{thm}
Comparing Theorem \ref{gf-Chen-D-f} with Theorem \ref{Sun-gene-thm}, we arrive at 
\begin{thm} \label{gf-refine} For $n\geq 0$, 
\begin{equation} \label{two-variable-chen-1}
G_{n+1}(x_1,x_2; y_1,y_2)=M_{n}(x_1y_1;x_2+y_2).
\end{equation}
\end{thm}

 Substituting $y_1=y_2=1$ into \eqref{two-variable-chen-1} yields Theorem \ref{chenrel}. 
Setting $x_1=x_2=x$ and $y_1=y_2=y$ in \eqref{two-variable-chen-1} leads to a result established by Chen and Pan \cite[Theorem 1.3 (1.14)]{Chen-Pan-2017} via combinatorial construction.

Setting $x_2=-y_2$ in \eqref{two-variable-chen-1} and using \eqref{eq-fourvar},  we obtain the following combinatorial assertions:  
\begin{cor} \label{cortt} {\rm (a)} For $n\geq 1$, the Catalan numbers $C_n$  count the number of plane trees with $2n+1$ edges and $n+1$ old leaves.

{\rm (b)} Let $P_{ey}(n,m)$ {\rm (}resp. $P_{oy}(n,m)${\rm )} denote the number of plane trees with $n$ edges, $m$ old leaves and an even (resp. odd) number of young leaves. For $n\geq 1$ and $\ 1\le m\le {\left\lfloor \frac{n}{2} \right\rfloor}$, 
\begin{align*}
{P}_{ey}(n,m)-{P}_{oy}(n,m) &=0.
\end{align*}
\end{cor}

\begin{remark}
   Corollary \ref{cortt} (a) can be proved by adding an old leaf on each vertex of plane trees with $n$ edges.
  It would be interesting to find a bijective proof of Corollary \ref{cortt} (b).
\end{remark}

Observe that 
$$
\begin{aligned}
M_n(x_1y_1;x_2+y_2) 
&= \sum_{k\geq 0}\binom{n}{2k} C_k (x_1y_1)^{k+1}(x_2+y_2)^{n-2k}\\[5pt]
&=\sum_{k\geq 0}\binom{n}{2k} C_k (x_1y_1)^{k+1} \sum_{i\geq 0} \binom{n-2k}{i}x_2^{n-2k-i}y_2^{i}\\[5pt]
&=\sum_{i\geq 0}\sum_{k\geq 0} \binom{n}{i}\binom{n-i}{2k}  C_k x_2^{n-2k-i}y_2^{i} (x_1y_1)^{k+1}\\[5pt]
&=\sum_{i=0}^{n}\binom{n}{i}y_2^i M_{n-i}(x_1y_1;x_2).
\end{aligned}
$$
Hence, we derive from \eqref{two-variable-chen-1} the following identity: 

\begin{cor}\label{Lin-Kim-g} For $n\geq 0$,  
   $$G_{n+1}(x_1,x_2;y_1,y_2) =  \sum_{i=0}^{n}\binom{n}{i}y_2^i M_{n-i}(x_1y_1;x_2).$$
\end{cor}
Setting $x_1=x_2=x$ and $y_1=y_2=1$ in Corollary \ref{Lin-Kim-g} results in the following identity established by Lin and Kim \cite{Lin-Kim-2022}:  
    \begin{equation}
    N_{n+1}(x)=\sum_{k=0}^{n}\binom{n}{k}M_{k}(x;x).
    \end{equation}
 When $x=1$, this polynomial reduces to the Euler transformation \eqref{rmotcat} (b) of the Motzkin numbers and the Catalan numbers.
 
By setting $y_1=1$ and $y_2=0$ in \eqref{two-variable-chen-1}, we obtain the following combinatorial interpretation of the Motzkin polynomials due to Donaghey\cite{Donaghey-1977}: 
\begin{equation}\label{com-Motzkin}
M_{n}(x_1;x_2)=G_{n+1}(x_1,x_2;1,0).
\end{equation}
More precisely, let $\mathcal{T}_n$ denote the set of plane trees with $n$ edges without  young interior vertices. According to \eqref{com-Motzkin}, we have 
\begin{equation}
    M_n(u;v) = \sum_{T\in \mathcal{T}_{n+1}} u^{ {\rm oleaf}(T)}v^{{\rm yleaf}(T)}.
\end{equation}
 This kind of plane trees was called by Donaghey\cite{Donaghey-1977} as the \textit{tip-augmented} plane trees. In other words, every interior vertex of a \textit{tip-augmented} plane tree is the parent of an old leaf.  For example,  only the last two plane trees in Fig. \ref{example-plane trees} are the tip-augmented plane trees, and so 
  \[M_2(u;v)=u^2+uv^2.\]
Furthermore, Fig. \ref{taptfoure} shows four tip-augmented plane trees with four edges, where the old leaves are labeled by $u$ and the young leaves are labeled by $v$. Therefore, we have 
\[M_3(u;v)=uv^3+3u^2v.\] 

  \begin{figure}[H]
    \centering
    \begin{tikzpicture}[scale=0.5]
        \draw[fill](0,0) circle(2pt)--(-1.5,-1)circle(2pt)node[below]{$u$};
        \draw[fill](0,0)--(-0.5,-1)circle(2pt)node[below]{$v$};
        \draw[fill](0,0)--(0.5,-1)circle(2pt)node[below]{$v$};
        \draw[fill](0,0)--(1.5,-1)circle(2pt)node[below]{$v$}; 
    \end{tikzpicture}\qquad
    \begin{tikzpicture}[scale=0.5]
        \draw[fill](0,0) circle(2pt)--(-1,-1)circle(2pt)node[below]{$u$};
        \draw[fill](0,0)--(0,-1)circle(2pt)--(0,-2)circle(2pt)node[below]{$u$};
        \draw[fill](0,0)--(1,-1)circle(2pt)node[below]{$v$}; 
\end{tikzpicture}\qquad
\begin{tikzpicture}[scale=0.5]
        \draw[fill](0,0) circle(2pt)--(-1,-1)circle(2pt)node[below]{$u$};
        \draw[fill](0,0)--(0,-1)circle(2pt)node[below]{$v$};
        \draw[fill](0,0)--(1,-1)circle(2pt)--(1,-2)circle(2pt)node[below]{$u$}; 
\end{tikzpicture}\qquad
     \begin{tikzpicture}[scale=0.5]
     \draw[fill](0,0) circle(2pt)--(-1,-1)circle(2pt)node[below]{$u$};
    \draw[fill](1,-1)--(0.5,-2)circle(2pt)node[below]{$u$};
    \draw[fill](0,0)--(1,-1)circle(2pt)--(1.5,-2)circle(2pt)node[below]{$v$}; 
    \end{tikzpicture}
    \caption{$4$-edge tip-augmented plane trees.} \label{taptfoure}
\end{figure}

To achieve further refinement of Coker's formula based on the polynomial   $G_n(x_{11},\break x_{12},x_2;y_{11},y_{12},y_2)$, we consider a refinement of the Motzkin polynomials. This refinement classifies the old leaves of a tip-augmented plane tree into three categories and the young leaves into two categories. More precisely,  
\begin{itemize}
\item A leaf without any siblings is said as a {\it  singleton leaf}; 

\item A leaf with siblings is considered  as an {\it  elder twin leaf}  if it is the leftmost child of its parent  and the second child of its parent is also a leaf; 

\item A leaf with siblings is considered  as an {\it  elder non-twin leaf} if it is the leftmost child of its parent  and the second child of its parent is not a leaf; 

\item  A leaf with siblings is considered  as a {\it second leaf} if it is the second child of its parent;

\item  A leaf with siblings is considered  as a {\it younger leaf} if it is neither the first child nor the second child of its parent; 
 
\end{itemize}
  Let $\text{\rm etleaf}(T)$, $\text{\rm entleaf}(T)$, $\text{\rm syleaf}(T)$ and $\text{\rm yerleaf}(T)$ denote the numbers of elder twin leaves,  elder non-twin leaves, second leaves and  younger leaves, respectively.

 \begin{defi}[Refinement of the Motzkin polynomials]
 For $n\geq 1$,
\begin{align*}
M_n(u_1,u_2,u_3;v_1,v_2)
&=\sum_{T\in\mathcal{T}_{n+1}} u_1^{{\rm sleaf}(T)}u_2^{\text {\rm etleaf}(T)}u_3^{\text {\rm entleaf}(T)}v_1^{\text{\rm yerleaf}(T)}v_2^{\text{\rm syleaf}(T)}
\end{align*}
with $M_0(u_1,u_2,u_3;v_1,v_2)=u_3$ by convention. 
\end{defi}
Fig. \ref{4-edge} shows further labeling of  four tip-augmented plane trees with four edges, where   the singleton leaves are labeled by $u_1$, the elder twin leaves are labeled by $u_2$,  the elder non-twin leaves are labeled by  $u_3$, the second leaves are labeled by  $v_2$, and  the younger leaves  are labeled by  $v_1$.  So 
\[{ M_3}(u_1,u_2,u_3; v_1,v_2)=u_2v_1^2v_2+u_1u_3v_1+u_1u_2v_2+u_2u_3v_2.\]

  \begin{figure}[H]
    \centering
    \begin{tikzpicture}[scale=0.5]
        \draw[fill](0,0) circle(2pt)--(-1.5,-1)circle(2pt)node[below]{$u_2$};
        \draw[fill](0,0)--(-0.5,-1)circle(2pt)node[below]{$v_2$};
        \draw[fill](0,0)--(0.5,-1)circle(2pt)node[below]{$v_1$};
        \draw[fill](0,0)--(1.5,-1)circle(2pt)node[below]{$v_1$}; 
    \end{tikzpicture}\qquad
    \begin{tikzpicture}[scale=0.5]
        \draw[fill](0,0) circle(2pt)--(-1,-1)circle(2pt)node[below]{$u_3$};
        \draw[fill](0,0)--(0,-1)circle(2pt)--(0,-2)circle(2pt)node[below]{$u_1$};
        \draw[fill](0,0)--(1,-1)circle(2pt)node[below]{$v_1$}; 
\end{tikzpicture}\qquad
\begin{tikzpicture}[scale=0.5]
        \draw[fill](0,0) circle(2pt)--(-1,-1)circle(2pt)node[below]{$u_2$};
        \draw[fill](0,0)--(0,-1)circle(2pt)node[below]{$v_2$};
        \draw[fill](0,0)--(1,-1)circle(2pt)--(1,-2)circle(2pt)node[below]{$u_1$}; 
\end{tikzpicture}\qquad
     \begin{tikzpicture}[scale=0.5]
     \draw[fill](0,0) circle(2pt)--(-1,-1)circle(2pt)node[below]{$u_3$};
    \draw[fill](1,-1)--(0.5,-2)circle(2pt)node[below]{$u_2$};
    \draw[fill](0,0)--(1,-1)circle(2pt)--(1.5,-2)circle(2pt)node[below]{$v_2$}; 
    \end{tikzpicture}
    \caption{$4$-edge tip-augmented plane trees.}\label{4-edge}
\end{figure}

\begin{thm}\label{thm-gf-overlineM} We have 
$$
\begin{aligned}
    &\quad\sum_{n\geq 0} M_n(u_1,u_2,u_3;v_1,v_2) q^{n}\\[5pt] 
    &= \frac{1-v_1q+(u_3-u_1)q^2}{2q^2}\\[5pt]
    &\quad-\frac{\sqrt{1-2v_1q+(v_1^2-2(u_1+u_3))q^2+ 2(u_1v_1+u_3v_1-2u_2v_2)q^3 +(u_3-u_1)^2q^4}}{2q^2}.
\end{aligned}
$$
\end{thm}
Setting $u_1=u_2=u_3=u$ and $v_1=v_2=v$ in Theorem \ref{thm-gf-overlineM} results in the generating function \eqref{Sun-gene}  established by Sun.  When $u_2=u_3=v_1=v_2=1$, we could recover the  generating function of the number of { (sharp)} peaks in Motzkin paths given by Brennan and Mavhungu \cite{Brennan-Mavhungu-2010} with the help of the bijection between the set of  tip-augmented plane trees and the set of elevated Motzkin paths, see \cite{Sulanke-2001}. 

By applying Theorem \ref{thm-gf-overlineM}, we find that  for $n\geq 2$, 
    $$M_n(u_1,u_2,u_3;v_1,v_2)= {M}_n(u_3,u_2,u_1;v_1,v_2),$$
which leads to the following combinatorial assertion: 
\begin{prop}\label{sym-tapt}
When $n \ge 3$, the number of tip-augmented plane trees with $n$ edges, $i$ singleton leaves, $j$ elder twin leaves, $k$ elder non-twin leaves, $r$ younger leaves and $s$ second leaves equals the number of tip-augmented plane trees with $n$ edges, $k$ singleton leaves, $j$ elder twin leaves, $i$ elder non-twin leaves, $r$ younger leaves and $s$ second leaves.
\end{prop}

{\remark It would be interesting to give a bijective proof of Proposition \ref{sym-tapt}. }

By comparing Theorem \ref{thm-gf-overlineM}  with Theorem \ref{thm-gf-abxy}, we derive the following refinement of Coker's formula.

\begin{thm} For $n\geq 1$,
    \[G_{n+1}(x_{11},x_{12},x_2;y_{11},y_{12},y_2)={M}_n(u_1,u_2,u_3;v_1,v_2)\]
 with     $u_2v_2=x_{12}y_{12}x_2+x_{11}y_{11}y_2$, $v_1=x_2+y_2$, $u_1=x_{11}y_{11}$ and $u_3=x_{12}y_{12}$.
\end{thm}
By setting $x_{11}=x_{12}=x_1$, $y_{11}=y_{12}=y_1$,  $u_1=u_2=u_3=u$ and $v_1=v_2=v$, we retrieve Theorem \ref{gf-refine}. 

As stated at the beginning of this paper, the primary objective of this paper is to introduce a grammatical approach to the study of the Narayana polynomials and the Motzkin polynomials. This technique, which employs context-free grammars to explore combinatorial polynomials, was pioneered by Chen \cite{Chen-1993}. 

In this paper, we  first establish the grammars for the  polynomial  $G_n(x_{11},x_{12},x_2;y_{11}, \break y_{12},
y_2)$ and the polynomial $M_n(u_1,u_2,u_3;v_1,v_2)$    (see Theorem \ref{thm-SOY-grammar} and Theorem \ref{thm-Refi-Motz}). We then provide  grammatical derivations for Theorem \ref{thm-gf-abxy} and Theorem \ref{thm-gf-overlineM}. As will be demonstrated later, the derivations of  the generating functions for $G_n(x_{11},x_{12},x_2;y_{11},y_{12},
5\break y_2)$ and $M_n(u_1,u_2,u_3;v_1,v_2)$  become quite simple once their grammars are established. In fact,  relying merely on grammars, we could deduce various generating functions and identities
without minding the recurrence relations or differential equations.
Additionally,  the grammatical approach is also applied to devise bijections \cite{Chen-Fu-2023, Chen-Fu-2024}, prove the $\gamma$-positivity of combinatorial polynomials  \cite{Chen-Fu-2022, Chen-Fu-Yan-2024, Lin-Ma-Zhang-2021, Ma-Ma-Yeh-2019}  and the stability of multivariable  combinatorial polynomials \cite{Chen-Hao-Yang-2021, Yang-Zhang-2024}.  In this regard, we anticipate further applications of these two grammars.   

It is worth noting that the grammar  for $G_n(x_{11},x_{12},x_2;y_{11}, y_{12}, y_2)$  presented in \eqref{SOY-grammar} and the grammar for $M_n(u_1,u_2,u_3;v_1,v_2)$ provided in  \eqref{Refi-Motz} can be specialized to yield    the grammars for the polynomial $G_n(x_{1},x_2;y_{1},
y_2)$ introduced by   Chen, Deutsch, and Elizalde, as well as the Motzkin polynomial $M_n(u;v)$ described in  Theorem \ref{thm-gramma-Chen-E-D} and Theorem \ref{gramar-Motzkin-thm}, respectively. Concerning the definition of grammars and backgrounds, please see Section 2.  The grammar for the Narayana polynomials has been established by Ma, Ma and Yeh \cite{Ma-Ma-Yeh-2019} and Yang and Zhang \cite{Yang-Zhang-2024}.  
  
\begin{thm}\label{thm-gramma-Chen-E-D}
Let $D$ be the formal derivative with respect to the grammar
\begin{equation}\label{gramma-Chen-E-D}
    G=\{x_1y_1\rightarrow 2tx_1y_1(x_2+y_2), x_2 \rightarrow 2tx_1y_1, y_2\rightarrow 2tx_1y_1, t\rightarrow t^2(x_2+y_2)\}. 
\end{equation}
For $n\geq 0$, 
\[D^{n}(y_2)=(n+1)!t^n{G}_{n}(x_1,x_2; y_1,y_2).\]
\end{thm}

\begin{thm}\label{gramar-Motzkin-thm}
Let $D$ be the formal derivative with respect to the grammar
\begin{equation} \label{gramar-Motzkin}
    M=\{t\rightarrow t^2v,~u\rightarrow 2tuv,~v\rightarrow 4tu\}. 
\end{equation}
For $n\geq 1$,  
\[D^{n}\bigg(\frac{v}{2}\bigg)=(n+1)!t^n{M}_{n-1}(u;v).\]
\end{thm}

 To conclude the  introduction, let us briefly discuss the real-rootedness  of the Motzkin polynomials and those defined by Chen, Deutsch and Elizalde \cite{Chen-Deutsch-Elizalde-2006}. This can be viewed as another application of Coker's formula. 

 Define 
 \begin{align}
 M_n(x):=\sum_{T \in \mathcal{T}_{n+1}}x^{{\rm oleaf}(T)} \quad \text{and} \quad 
 G_n(x):=\sum_{T \in \mathcal{P}_{n}}x^{{\rm oleaf}(T)},
 \end{align}
 where $\mathcal{T}_n$ denotes the set of  tip-augmented  plane trees with $n$ edges and $\mathcal{P}_n$ denotes the set of    plane trees with $n$ edges. 
We have the following consequence: 
 \begin{thm}\label{real-roots} For $n\geq 1$, the polynomial $G_n(x)$ and  $M_n(x)$ have
     only  non-positive real zeros.
\end{thm}
 
 Theorem  \ref{real-roots} can be derived  by applying the fact that the Narayana polynomial $N_n(x)$ has only non-positive real zeros, Coker's formula   \eqref{two-variable},    the relation \eqref{eqchenrel}  and the following observation stated by Petersen \cite[Observation 4.2]{Petersen-2015}.

\begin{prop}{\rm \cite[Observation 4.2]{Petersen-2015}}\label{thm-6-25}
    If a polynomial $f(x)$ is symmetric,
then $$f(x)=\sum_{k=0}^{\lfloor \frac{n}{2} \rfloor } \gamma_k x^k (1+x)^{n-2k}$$ 
has only negative (non-positive) real zeros if and only if $$\gamma(f;x)=\sum_{k=0}^{\lfloor \frac{n}{2} \rfloor }  \gamma_k x^k$$ has only negative (non-positive) real zeros.
\end{prop}

This paper is organized as follows. Section 2 aims to illustrate the idea of the grammatical approach through a grammatical calculus  for the Motzkin polynomials. In Section 3, we provide a grammatical derivation of Theorem \ref{thm-gf-abxy}. This involves establishing the grammar for the polynomial  $G_{n}(x_{11},x_{12},x_2;y_{11},y_{12},y_2)$ (see Theorem \ref{thm-SOY-grammar}). Section 4 focuses on the grammatical derivation of Theorem  \ref{thm-gf-overlineM}. The grammar for the polynomial $ M_n(u_1,u_2,u_3;v_1,v_2)$ is detailed in  Theorem \ref{thm-Refi-Motz}.

\section{ A grammatical calculus for the Motzkin polynomials}

To demonstrate how  the grammatical approach works, we first prove Theorem \ref{gramar-Motzkin-thm} by using the grammatical labeling and then  we derive the generating function for $M_n(u;v)$ (Theorem \ref{Sun-gene-thm}) solely based on the grammar  in Theorem \ref{gramar-Motzkin-thm}.

A context-free grammar $G$ over a set $V=\{x,y,z,\ldots\}$ of variables is a set of substitution rules replacing a variable in $V$ by a Laurent polynomial of variables in $V$. For a context-free grammar $G$ over $V$, the formal derivative $D$ with respect to $G$ is defined as a linear operator acting on Laurent polynomials with variables in $V$ such that each substitution rule is treated as the common differential rule that satisfies  the following relations:
\begin{align} \label{gramma-add-rela}
D(u+v)&=D(u)+D(v),\\[5pt]
D(uv)&=D(u)v+uD(v). \label{gramma-mult-rela}
\end{align}
Hence, it obeys the Leibniz's rule
\[D^{n}(uv)=\sum_{k=0}^n\binom{n}{k}D^k(u)D^{n-k}(v).
\]
For a constant $c$, we have $D(c)=0$.

A formal derivative $D$ with respect to $G$ is also associated with an exponential generating function. For a  Laurent polynomial $w$ of variables in $V$, let
\begin{equation}
    {\rm Gen}(w;q)=\sum_{n\geq 0}D^n(w)\frac{q^n}{n!}.
\end{equation}
Then, by \eqref{gramma-add-rela} and \eqref{gramma-mult-rela}, we derive that
\begin{equation}\label{gramma-add}
{\rm Gen}(u+v;q) = {\rm Gen}(u;q)+{\rm Gen}(v;q),
\end{equation}
\begin{equation}\label{gramma-multiple}
{\rm Gen}(uv;q) = {\rm Gen}(u;q){\rm Gen}(v;q).
\end{equation}
For more information on the grammatical calculus, we refer to Chen \cite{Chen-1993} and Chen and Fu \cite{Chen-Fu-2017, Chen-Fu-2022}.

Next, we apply the grammatical labeling to  prove Theorem \ref{gramar-Motzkin-thm}. The
notion of a grammatical labeling was introduced by Chen and Fu \cite{Chen-Fu-2017}.   Here we need to consider  labeled  tip-augmented plane trees. Recall that a tip-augmented plane tree is a plane tree without young interior vertices. A labeled plane tree with  $n$ edges is a  plane tree where each vertex is uniquely labeled with a number from the set 
$[n+1]=\{1,2,\dots,n+1\}$.
 Let ${\mathcal{LT}_{n}}$ denote the set of all labeled tip-augmented plane trees with $n$ edges. 

For $T\in{\mathcal{LT}_{n}}$, we define the grammatical labeling of $T$ as follows: 
\begin{itemize}
\item If the vertex $\mathbf{j}$  of $T$ is an old leaf, then label $\mathbf{j}$ by $u$;

\item If the vertex  $\mathbf{j}$ of $T$  is a young leaf, then label  $\mathbf{j}$ by $v$;

\item All of  edges of $T$  are labeled by $t$.
\end{itemize}

The weight of $T$ is defined to be the product of all the labels, that is, 
\[{\rm wt}(T)=u^{ {\rm oleaf}(T)} v^{ {\rm yleaf}(T)}t^{{\rm edge}(T)}.\]

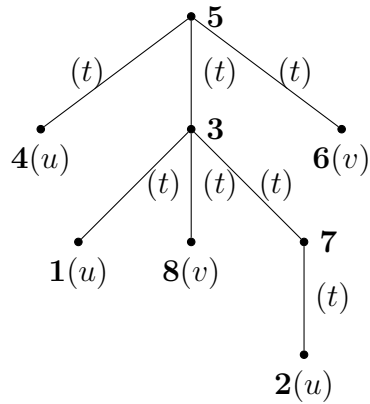
\begin{figure}[H]
    \centering
        \begin{tikzpicture}
[vertex/.style={shape=circle, draw, inner sep=1pt, fill=black},
sibling distance=2cm,level distance=1.5cm]
\node[vertex,label=0:{$\mathbf{5}$}]{}[grow=down]
  child {node [vertex, label=-90:{$\mathbf{4}(u)$}]{}edge from parent node[left]{$(t)$}}
   child {node [vertex, label=0:{$\mathbf{3}$}]{}[sibling distance=1.5cm]
   child {node [vertex, label=-90:{$\mathbf{1}(u)$}]{}edge from parent node[right]{$(t)$}}
    child {node [vertex, label=-90:{$\mathbf{8}(v)$}]{}edge from parent node[right]{$(t)$}}
   child {node [vertex, label=0:{$\mathbf{7}$}]{}
    child {node [vertex, label=270:{$\mathbf{2}(u)$}]{}edge from parent node[right]{$(t)$}
   }edge from parent node[right]{$(t)$}}edge from parent node[right]{$(t)$}}
    child {node [vertex, label=270:{$\mathbf{6}(v)$}]{}edge from parent node[right]{$(t)$}};
\end{tikzpicture}
       \caption{A    tip-augmented plane tree in ${\mathcal{LT}}_7$ with the $(u,v;t)$-labeling.}
    \label{tiplizi}
     \end{figure}

For example, Fig. \ref{tiplizi} shows the grammatical labeling of a labeled  tip-augmented plane tree $T\in {\mathcal{LT}_7
}$ whose weight of $T$ is ${\rm wt}(T)=t^7u^3v^2$. The
labels are shown in parentheses. We refer to this labeling scheme of    tip-augmented plane trees as 
the $(u,v;t)$-labeling. 

From the definition of the $(u,v;t)$-labeling, we see that the polynomial $M_{n}(u;v)$ can be interpreted as:
\begin{equation}
    (n+2)!t^{n+1}M_{n}(u;v)=\sum_{T\in \mathcal{LT}_{n+1}}{\rm wt}(T).
\end{equation}
Therefore, the proof of Theorem \ref{gramar-Motzkin-thm} is equivalent to the proof of the following assertion: For $n\geq 1$, 
\begin{equation}\label{gram-induct-motzkin}
    D^{n}\bigg(\frac{v}{2}\bigg) = \sum_{T\in \mathcal{LT}_{n}}{\rm wt}(T),
\end{equation}
     where $D$ is the formal derivative with respect to the grammar defined in \eqref{gramar-Motzkin}.

\begin{figure}[H]
    \centering
     \begin{tikzpicture}
[vertex/.style={shape=circle, draw, inner sep=1pt, fill=black},
sibling distance=1.5cm,level distance=1.5cm,]
\node[vertex, label=0:{$\mathbf{1}$}]{}[grow=down]
  child {node [vertex, label=0:{$\mathbf{2}(u)$}]{} edge from parent node[right]{$(t)$}};
\end{tikzpicture}\qquad
 \begin{tikzpicture}
[vertex/.style={shape=circle, draw, inner sep=1pt, fill=black},
sibling distance=1.5cm,level distance=1.5cm,]
\node[vertex, label=0:{$\mathbf{2}$}]{}[grow=down]
  child {node [vertex, label=0:{$\mathbf{1}(u)$}]{} edge from parent node[right]{$(t)$}};
\end{tikzpicture}
    \caption{Two tip-augmented plane trees in $\mathcal{LT}_1$ with the $(u,v;t)$-labeling. } \label{exgraMot}
\end{figure}
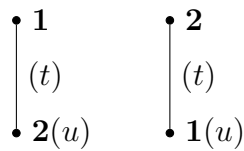

\noindent{\it Proof of Theorem \ref{gramar-Motzkin-thm}.} We proceed by induction on $n$. For $n=1$, there are two labeled tip-augmented plane trees. The $(u,v;t)$-labelings of those two   trees are given in Fig. \ref{exgraMot}.  It is easy to check that the relation \eqref{gram-induct-motzkin} is valid for $n=1$. 
Assume that this assertion holds for $n-1$, that is, the relation \eqref{gram-induct-motzkin} is valid for $n-1$. To demonstrate that it also holds for $n$,   it suffices to show, by \eqref{gram-induct-motzkin}, that 
\begin{equation}\label{gram-inductaa1}
D\left(\sum_{T\in \mathcal{LT}_{n-1}}{\rm wt}(T)\right)=\sum_{T\in \mathcal{LT}_{n}}{\rm wt}(T),
\end{equation}
where $D$ is the formal derivative with respect to the grammar defined in \eqref{gramar-Motzkin}.

We will accomplish this by constructing all labeled tip-augmented plane trees $\widetilde{T}$ through the insertion of the vertex $\mathbf{n+1}$ into labeled tip-augmented plane trees  in $\mathcal{LT}_{n-1}$, ensuring that the change of weights between these two plane trees adheres to the substitution rules given in the grammar \eqref{gramar-Motzkin}.

Given a tip-augmented plane   tree $T\in \mathcal{LT}_{n-1}$, and let $\widetilde{T}$ be the  plane tree with $n+1$ vertices  obtained by performing a specific operation at the position $\diamond$ of $T$.  We consider the following three cases: 

\begin{itemize}
 \item[(1).]  The position $\diamond$  is  an old leaf $\mathbf{j}$ of $T$ labeled by  $u$.  Suppose that $\mathbf{i}$ is the parent of $\mathbf{j}$, see Fig. \ref{(I)}.  Add the vertex $\mathbf{n+1}$ as illustrated in Fig. \ref{(I-i)} to construct the plane tree $\widetilde{T}$. It is easy to check that the change of weights follows the substitution rule $u\rightarrow 2tuv$. Moreover, every interior vertex of the plane tree $\widetilde{T}$  has an old leaf, so the plane tree $\widetilde{T}$ is a tip-augmented plane tree. 

 \begin{figure}[H]
    \centering
     \begin{tikzpicture}
[vertex/.style={shape=circle, draw, inner sep=1pt, fill=black},
subtree/.style={shape=ellipse, draw,minimum width=1.5cm, minimum height=.5cm},
sibling distance=1.5cm,level distance=0.8cm,
leaf/.style={label={[name=#1]below:$#1$}},auto]
\node{}[grow=down]
  child {node [vertex, label=0:{$\mathbf{i}$}]{}[level distance=10mm]
   child {node [diamond, draw, solid, inner sep=1.3pt, label=270:{$\mathbf{j}(u)$}]{} edge from parent [solid]node[left]{$(t)$}
}
     child {node [subtree, solid, rotate=126]{}  edge from parent [solid]}
edge from parent [dashed]};
\end{tikzpicture}
    \caption{The position
$\diamond$ is an old leaf of $T$.}
    \label{(I)}
\end{figure}
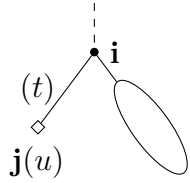

\begin{figure}[H]
    \centering
   \begin{subfigure}[t]{0.3\textwidth}
   \centering 
        \begin{tikzpicture}
[vertex/.style={shape=circle, draw, inner sep=1pt, fill=black},
subtree/.style={shape=ellipse, draw,minimum width=1.5cm, minimum height=.5cm},
sibling distance=1.5cm,level distance=0.8cm,
leaf/.style={label={[name=#1]below:$#1$}},auto]
\node{}[grow=down]
  child {node [vertex, label=0:{$\mathbf{i}$}]{}[level distance=15mm]
   child {node [vertex, label=270:{$\uline{\mathbf{n\!+\!1}}(u)$}]{}edge from parent [solid] node[left,xshift=-1mm]{$(t)$} }
   child {node [vertex, label=270:{$\mathbf{j}(v)$}]{}  edge from parent [solid]node[left,xshift=1mm]{$(t)$}}
   child {node [subtree,solid, rotate=128]{}  edge from parent [solid]}edge from parent [dashed]
};
\end{tikzpicture}
        \caption{$\mathbf{n+1}$ is an old leaf of $\widetilde{T}$ with at least one sibling, and the sibling immediately following $\mathbf{n+1}$ is a young leaf.}
        \label{(I-1)}
    \end{subfigure}\qquad 
    \begin{subfigure}[t]{0.3\textwidth}
    \centering 
        \begin{tikzpicture}
[vertex/.style={shape=circle, draw, inner sep=1pt, fill=black},
subtree/.style={shape=ellipse, draw,minimum width=1.5cm, minimum height=.5cm},
sibling distance=1.5cm,level distance=0.8cm,
leaf/.style={label={[name=#1]below:$#1$}},auto]
\node{}[grow=down]
  child {node [vertex, label=0:{$\uline{\mathbf{n\!+\!1}}$}]{}[level distance=15mm]
   child {node [vertex, label=270:{$\mathbf{i}(u)$}]{}edge from parent [solid]node[left,xshift=-1mm]{$(t)$} }
   child {node [vertex, label=270:{$\mathbf{j}(v)$}]{}  edge from parent [solid]node[left,xshift=1mm]{$(t)$}}
   child {node [subtree,solid, rotate=128]{}  edge from parent [solid]}edge from parent [dashed]
};
\end{tikzpicture}
        \caption{$\mathbf{n+1}$ is an interior vertex of  $\widetilde{T}$ with  at least two children and its second child  is a young leaf.}
        \label{(I-2)}
    \end{subfigure}
    
     \caption{$u\rightarrow 2tuv$.}
     \label{(I-i)}
\end{figure}
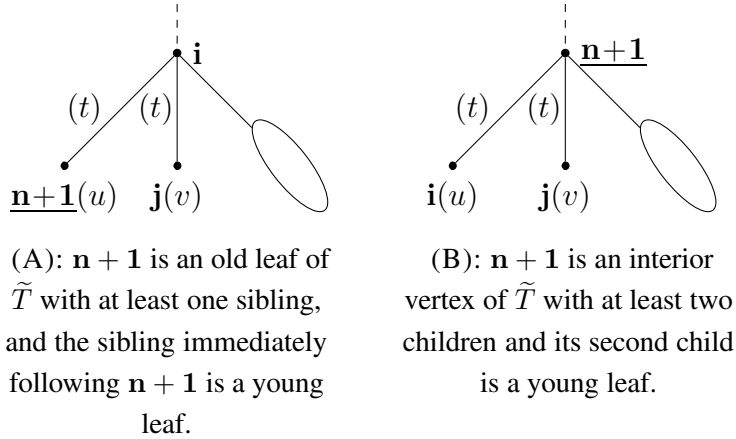

\item[(2).]  The position $\diamond$  is  a young leaf $\mathbf{j}$ of $T$ labeled by  $v$, see Fig. \ref{(II)}.   It should be noted that the forest $A$  formed  by the subtrees $\mathbf{k_{2}}, \ldots, \mathbf{k_{t-1}}$  and the forest $B$ formed by the subtrees $\mathbf{k_{t+1}}, \ldots, \mathbf{k_{q}}$ may be empty. We perform the steps shown in Fig. \ref{(II-i)} to create four kinds of tip-augmented plane trees. It is easy to check that the change of weights  follows the substitution rule $v\rightarrow 4tu$.

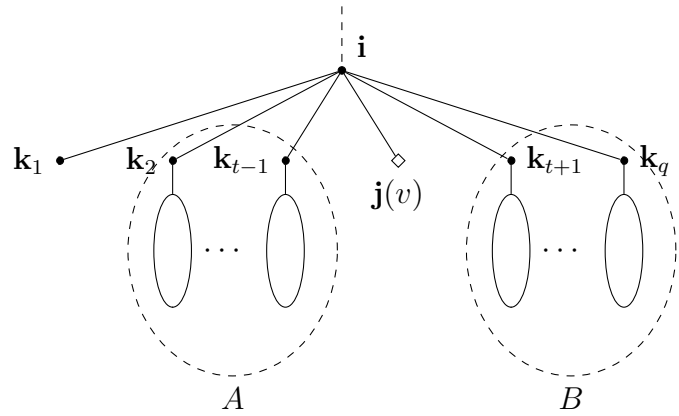
\begin{figure}[H]
    \centering
      \begin{tikzpicture}
[vertex/.style={shape=circle, draw, inner sep=1pt, fill=black},
subtree/.style={shape=ellipse, draw,minimum width=1.5cm, minimum height=.5cm},
every fit/.style={ellipse,draw,inner sep=-2pt},
sibling distance=1.5cm,level distance=1cm,
leaf/.style={label={[name=#1]below:$ $}}]

\node{}[grow=down]
  child {node [vertex, label=30:{$\mathbf{i}$}]{}[level distance=12mm]
  child {node [vertex,label=180:{$\mathbf{k}_1$}]  {} edge from parent [solid]}
  child {node [vertex, label=180:{$\mathbf{k}_2$}] (a's parent) {}  [level distance=12mm]
   child {node [subtree, rotate=90,leaf=a,label=-90:{$\cdots$}]{} [level distance=9mm] child {node [ rotate=70,leaf=c]{} edge from parent [solid]}edge from parent [solid]}edge from parent [solid]}
   child {node [vertex, label=180:{$\mathbf{k}_{t-1}$}] (b's parent) {}  [level distance=12mm]
   child {node [subtree, rotate=90,leaf=b]{} edge from parent [solid]}edge from parent [solid]}   
   child {node [diamond, draw, solid, inner sep=1.3pt,label=270:{$\mathbf{j}(v)$}]{} edge from parent [solid]} 
   child {node [vertex, label=0:{$\mathbf{k}_{t+1}$}] (d's parent) {}  [level distance=12mm]
   child {node [subtree, rotate=90,leaf=d,label=-90:{$\cdots$}]{}[level distance=9mm] child {node [ rotate=70,leaf=f]{}edge from parent [solid] }edge from parent [solid] }edge from parent [solid]}
   child {node [vertex, label=0:{$\mathbf{k}_q$}] (e's parent) {}  [level distance=12mm]
   child {node [subtree, rotate=90,leaf=e]{} edge from parent [solid]}edge from parent [solid]}
   edge from parent [dashed]};
   \node[draw,dashed,xshift=-2mm,fit=(a's parent)(a) (b) (c)  ,label=below:$A$] {}; 
     \node [dashed,xshift=-2mm,fit=(e) (f) (d's parent),label=below:$B$] {};
\end{tikzpicture}
        \caption{The position $\diamond$ is a young leaf of $T$.}
        \label{(II)}
\end{figure}

\begin{figure}[H]
    \centering
    \begin{subfigure}[t]{0.4\textwidth}
    \centering
        \begin{tikzpicture}
[vertex/.style={shape=circle, draw, inner sep=1pt, fill=black},
subtree/.style={shape=ellipse, draw,minimum width=1.5cm, minimum height=.5cm},
every fit/.style={ellipse,draw,inner sep=-2pt},
sibling distance=1.5cm,level distance=1cm,
leaf/.style={label={[name=#1]below:$ $}},auto]

\node{}[grow=down]
  child {node [vertex, label=30:{$\mathbf{i}$}]{}[level distance=12mm]
  child {node [vertex,label=180:{$\mathbf{k}_1$}]{}edge from parent [solid]}
  child {node [subtree, rotate=60,label=190:{$A$}]{}edge from parent [solid]}  
   child {node [vertex, label=0:{$\mathbf{j}$}]{} [level distance=12mm]
   child {node [vertex, label=270:{$\uline{\mathbf{n\!+\!1}}(u)$}]{}edge from parent [solid]node[right]{$(t)$}} edge from parent[solid]} 
   child {node [subtree, rotate=-30,label=-60:{$B$}]{}edge from parent [solid]} 
   edge from parent [dashed]};
\end{tikzpicture}
        \caption{$\mathbf{n+1}$ is an old leaf of $\widetilde{T}$ with no siblings.}
        \label{(II-1)}
    \end{subfigure}\qquad
  \begin{subfigure}[t]{0.4\textwidth}
    \centering        \begin{tikzpicture}
[vertex/.style={shape=circle, draw, inner sep=1pt, fill=black},
subtree/.style={shape=ellipse, draw,minimum width=1.5cm, minimum height=.5cm},
every fit/.style={ellipse,draw,inner sep=-2pt},
sibling distance=1.5cm,level distance=1cm,
leaf/.style={label={[name=#1]below:$ $}},auto]

\node{}[grow=down]
  child {node [vertex, label=30:{$\mathbf{i}$}]{}[level distance=12mm]
  child {node [vertex,label=180:{$\mathbf{k}_1$}]{}edge from parent [solid]}
  child {node [subtree, rotate=60,label=190:{$A$}]{}edge from parent [solid]}  
   child {node [vertex, label=180:{$\uline{\mathbf{n\!+\!1}}$}]{} [level distance=12mm]
   child {node [vertex, label=270:{$\mathbf{j}(u)$}]{}edge from parent [solid]node[right]{$(t)$}} edge from parent[solid]} 
   child {node [subtree, rotate=-30,label=-60:{$B$}]{}edge from parent [solid]} 
   edge from parent [dashed]};
\end{tikzpicture}
        \caption{$\mathbf{n+1}$ is an interior vertex   of $\widetilde{T}$ with  only one child.}
        \label{(II-2)}
    \end{subfigure}

\begin{subfigure}[t]{0.4\textwidth}
   \centering
        \begin{tikzpicture}
[vertex/.style={shape=circle, draw, inner sep=1pt, fill=black},
subtree/.style={shape=ellipse, draw,minimum width=1.5cm, minimum height=.5cm},
sibling distance=1.5cm,level distance=1.2cm,
leaf/.style={label={[name=#1]below:$#1$}},auto]
\node{}[grow=down]
  child {node [vertex, label=10:{$\mathbf{j}$}]{}[level distance=12mm]
  child {node [vertex, label=270:{$\uline{\mathbf{n\!+\!1}}(u)$}]{}edge from parent[solid]node[left,yshift=1mm]{$(t)$}}
  child {node [vertex, label=180:{$\mathbf{i}$}]{}
  child {node[vertex,label=270:{$\mathbf{k}_1$}] {} edge from parent [solid]}
   child {node [subtree,  draw, dashed , rotate=120,label=205:{$A$}]{} edge from parent [solid]}
   edge from parent [solid] }
   child {node [subtree,  draw, dashed ,rotate=140,label=170:{$B$}]{} edge from parent [solid]}edge from parent [dashed]
};
\end{tikzpicture}
        \caption{$\mathbf{n+1}$ is an old leaf of $\widetilde{T}$ with  at least one sibling,  and the sibling immediately following $\mathbf{n+1}$ is not a leaf.  }
        \label{(II-3)}
    \end{subfigure}\qquad
    \begin{subfigure}[t]{0.4\textwidth}
   \centering
        \begin{tikzpicture}
[vertex/.style={shape=circle, draw, inner sep=1pt, fill=black},
subtree/.style={shape=ellipse, draw,minimum width=1.5cm, minimum height=.5cm},
sibling distance=1.5cm,level distance=1.2cm,
leaf/.style={label={[name=#1]below:$#1$}},auto]
\node{}[grow=down]
  child {node [vertex, label=10:{$\uline{\mathbf{n\!+\!1}}$}]{}[level distance=12mm]
  child {node [vertex, label=270:{$\mathbf{j}(u)$}]{}edge from parent[solid]node[left,yshift=1mm]{$(t)$}}
  child {node [vertex, label=180:{$\mathbf{i}$}]{}
  child {node[vertex,label=270:{$\mathbf{k}_1$}] {} edge from parent [solid]}
   child {node [subtree,  draw, dashed , rotate=120,label=205:{$A$}]{} edge from parent [solid]}
   edge from parent [solid] }
   child {node [subtree,  draw, dashed ,rotate=140,label=170:{$B$}]{} edge from parent [solid]}edge from parent [dashed]
};
\end{tikzpicture}
        \caption{$\mathbf{n+1}$ is an interior vertex  of $\widetilde{T}$ with at least two children and its second child  is not a  leaf.}
        \label{(II-4)}
    \end{subfigure}
  
     \caption{$v\rightarrow 4tu$.}
     \label{(II-i)}
\end{figure}

\item[(3)] The position $\diamond$ is an edge $(\mathbf{i},\mathbf{j})$ in $T$ labeled by $t$, as shown in Fig. \ref{(III)}.  Insert the edge $(\mathbf{i},\mathbf{n+1})$ right immediately  after $(\mathbf{i},\mathbf{j})$ as illustrated in Fig. \ref{(III-i)} to get the tip-augmented plane tree  $\widetilde{T}$. The change of weights is accordance with  $t\rightarrow t^2v$. 

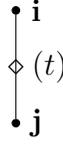
\begin{figure}[H]
    \centering
    \begin{tikzpicture}
 [vertex/.style={shape=circle, draw, inner sep=1pt, fill=black},
 subtree/.style={shape=ellipse, draw,minimum width=1.5cm, minimum height=.5cm},
 every fit/.style={ellipse,draw,inner sep=-2pt},
sibling distance=1.8cm,level distance=1.5cm,
 leaf/.style={label={[name=#1]below:$ $}},auto]
\node[vertex,label=0:{$\mathbf{i}$}]{}[grow=down]
  child {node [vertex, label=0:{$\mathbf{j}$}]{}[level distance=20mm]
 edge from parent [solid] node[xshift=-2.5mm] {$\diamond ~(t)$}};
 
 \end{tikzpicture}
   \caption{The position $\diamond$ is an  edge.}
  \label{(III)}
\end{figure}

\begin{figure}[H]
    \centering
    \begin{subfigure}[t]{0.4\textwidth}
    \centering
    \begin{tikzpicture}
 [vertex/.style={shape=circle, draw, inner sep=1pt, fill=black},
 subtree/.style={shape=ellipse, draw,minimum width=1.5cm, minimum height=.5cm},
 every fit/.style={ellipse,draw,inner sep=-2pt},
sibling distance=1.8cm,level distance=1.5cm,
 leaf/.style={label={[name=#1]below:$ $}},auto]
\node[vertex,label=0:{$\mathbf{i}$}]{}[grow=down]
  child {node [vertex, label=-90:{$\mathbf{j}$}]{}[level distance=20mm]
 edge from parent [solid] node[left] {$(t)$} }
 child{node[vertex, label=-90:{$\uline{\mathbf{n\!+\!1}}(u)$}]{} edge from parent [solid] node[right]{$(t)$}};
 
 \end{tikzpicture}
 \caption{$\mathbf{n+1}$ is a young leaf of $\widetilde{T}$.}
 \end{subfigure}
    \caption{$t\rightarrow t^2 v$.}
   
  \label{(III-i)}
\end{figure}

 \end{itemize}

From the above construction, it is clear that 
\begin{itemize}
     \item  $\mathbf{n+1}$ is an old   leaf of $\widetilde{T}$ in Fig. \ref{(I-1)}, Fig. \ref{(II-1)} and Fig. \ref{(II-3)} .
    \item   $\mathbf{n+1}$ is a young leaf  of $\widetilde{T}$ in  Fig. \ref{(III-i)}.
     \item  $\mathbf{n+1}$ is an interior vertex of  $\widetilde{T}$ in Fig. \ref{(I-2)}, Fig. \ref{(II-2)} and Fig. \ref{(II-4)}.
\end{itemize}
Given this observation,  we see that the construction is  reversible depending on the position of the vertex $\mathbf{n+1}$ in $\widetilde{T}$. 
Consequently, we could generate 
all labeled tip-augmented plane trees in $\mathcal{LT}_{n}$ from those in $\mathcal{LT}_{n-1}$ using the above construction. Moreover, the change of weights between   tip-augmented plane trees in $\mathcal{LT}_{n}$  and those in $\mathcal{LT}_{n-1}$   is consistent with the substitution rule given in the grammar \eqref{gramar-Motzkin}. Thus we show that 
\eqref{gram-induct-motzkin} is valid and so the assertion is also valid for $n$. This completes the proof.  \qed

Equipped with Theorem \ref{gramar-Motzkin-thm}, we are now in a position to  provide a grammatical derivation of the generating function for $M_n(u;v)$ as stated in Theorem \ref{Sun-gene-thm}. 

 \noindent{ \bf The grammatical  derivation for Theorem \ref{Sun-gene-thm}.} For the formal derivative $D$ with respect to the grammar $M$ in \eqref{gramar-Motzkin}, 
  we observe that 
\[D(t^{-1}) = -v \quad \text{and} \quad
    D^{2}(t^{-1}) = -D(v).\]
Combining this with Theorem  \ref{gramar-Motzkin-thm},  we deduce that for $n\geq 2$, 
$$D^{n}(t^{-1})=-2n!t^{n-1}M_{n-2}(u;v).$$
Consequently, 
\begin{align} \label{eqnaa}
    {\rm Gen}(t^{-1};q) &= \sum_{n\geq 0} D^{n}(t^{-1})\frac{q^n}{n!}\nonumber \\[5pt]
    &=t^{-1}-vq-\sum_{n\geq 2}D^{n-1}\left(v\right)\frac{q^n}{n!} \nonumber \\[5pt] 
     &=t^{-1}-vq-2\sum_{n\geq 2}n!t^{n-1}M_{n-2}(u;v)\frac{q^n}{n!} \nonumber \\[5pt]
     &=t^{-1}-vq-2tq^2\sum_{n\geq 0} M_{n}(u;v)t^{n}q^n.
\end{align}
Since  
\begin{equation}\label{gft-1}
{\rm Gen}(t^{-2};q)={\rm Gen}^2(t^{-1};q),
\end{equation}
it suffices to compute ${\rm Gen}(t^{-2};q)$. Observe that 
$$D(t^{-2}) = -2t^{-1}v, \quad
    D^{2}(t^{-2}) = 2v^2-8u, \quad \text{and} \quad D^{3}(t^{-2}) = 0.
$$
It follows that
\begin{align}\label{eqnaaddtt}
    {\rm Gen}(t^{-2};q) &= \sum_{n\geq 0} D^{n}(t^{-2})\frac{q^n}{n!}=t^{-2}-2t^{-1}vq+(v^2-4u)q^2. 
\end{align}
By setting  $q=0$ in \eqref{gft-1},  we infer from \eqref{eqnaa} and  \eqref{eqnaaddtt} that  
\begin{equation}\label{gft-aa}
{\rm Gen}(t^{-1};q)= \sqrt{t^{-2}-2t^{-1}vq+(v^2-4u)q^2}.  
\end{equation}
Substituting \eqref{gft-aa} into \eqref{eqnaa} and setting $t=1$, we arrive at the generating function for  $M_n(u;v)$ as stated in Theorem \ref{Sun-gene-thm}. \qed

\section{{Proof of Theorem \ref{thm-gf-abxy}}}

The main objective of this section is to  give a   proof of Theorem  \ref{thm-gf-abxy} by using the grammatical calculus. 
To this end, we need to establish the grammar for the polynomial $G_n({x}_{11},{{{x}}}_{12},x_2;{y}_{11},{{y}}_{12},y_2)$.  It turns out that  we should classify the edge of a plane tree. An edge is called a young edge if it is not an edge of an old leaf (including a singleton leaf and an elder leaf). Let ${\rm yedge}(T)$ denote the number of young edges of a plane tree $T$.  

For $n\geq 2$, we define 
\begin{equation}\label{defi-G}
\widetilde{G}_n(a,b,c,d,t)=\sum_{T\in \mathcal{P}_{n}}a^{{\rm sleaf}(T)} b^{{\rm eleaf}(T)}c^{{\rm yleaf}(T)} d^{{\rm yint}(T)}t^{{\rm yedge}(T)}
\end{equation}
with the convention that 
  $\widetilde{G}_0(a,b,c,d,t)=d$ and $\widetilde{G}_1(a,b,c,d,t)=b$. 

By definition, it is evident that  for $n\geq 2$, 
\[\widetilde{G}_n(x_{11}y_{11}t,x_{12}y_{12}t,x_2,y_2,t)=t^nG_n(x_{11},x_{12},x_{2};y_{11},y_{12},y_{2}). \]

In this section,   we first show that the following grammar
  \begin{equation}\label{SOY-grammar}
    G=\{a\rightarrow 3t(ad+bc),~b\rightarrow 3t(ad+bc),~c\rightarrow 2a,~ d\rightarrow 2b,~t\rightarrow t^2(c+d)\} \end{equation}
 can be used to generate the polynomial $\widetilde{G}_n(a,b,c,d,t)$. More precisely,
\begin{thm}\label{thm-SOY-grammar}
Let $D$ be the formal derivative with respect to the grammar defined in \eqref{SOY-grammar}. For  $n\geq 2$,  
\[D^{n}(d)=(n+1)!\widetilde{G}_{n}(a,b,c,d;t).\]
\end{thm}
The first few values of $D^{n}(d)$ are given below:
$$ 
\begin{aligned}
D^2(d)&=6(tad+tbc),\\[5pt]
D^3(d)&=24(tab+t^2acd+t^2ad^2+t^2bc^2 + t^2bcd),\\[5pt]
D^4(d)&= 
120(t^2a^2d+2t^2abd+t^3ad^3+t^3ac^2d+2t^3acd^2+2t^2abc\\[5pt]
&\quad +t^2b^2c+t^3bc^3+2t^3bc^2d+t^3bcd^2).
\end{aligned}
$$

By setting $a=b=x_1y_1t$, $c=x_2$, $d=y_2$ in the grammar $G$ given in \eqref{SOY-grammar}, we obtain the grammar $G$ presented in \eqref{gramma-Chen-E-D} for the 
 polynomial $G_n(x_1,x_{2},y_1,y_{2})$. Thus, we derive Theorem \ref{thm-gramma-Chen-E-D} from Theorem \ref{thm-SOY-grammar}.

With Theorem \ref{thm-SOY-grammar}  in hand,  we obtain the following generating function for $\widetilde{G}_{n}(a,b,c,d;t)$. By setting $t=1$, $a=x_{11}y_{11},b=x_{12}y_{12}$, $c=x_2$ and $d=y_2$, it  can be specialized to obtain 
the generating function for $G_{n}(x_{11},x_{12},x_2;y_{11},y_{12},y_2)$   as stated in Theorem \ref{thm-gf-abxy}.

\begin{thm}\label{gf-gwilde}
$$
\begin{aligned}
&\sum_{n\geq 0} \widetilde{G}_{n}(a,b,c,d;t) q^n = \frac{t^{-1}+(d-c)q+(b-a)q^2}{2q}\\[5pt]
&\quad -\frac{\sqrt{ t^{-2}-2t^{-1}(c+d)q+((c+d)^2-2t^{-1}(a+b))q^2 + 2(a-b)(c-d)q^3+(a-b)^2q^4}}{2q}.
\end{aligned}
$$
\end{thm}

 \subsection{A grammatical labeling of $\widetilde{G}_{n}(a,b,c,d;t)$} 
 In this subsection, we aim to  show Theorem \ref{thm-SOY-grammar} by using the grammatical labeling. Similarly, we need to consider  labeled plane trees.  Let $\mathcal{LP}_{n}$ denote the set of all labeled plane trees with $n$ edges. 

The grammatical labeling of  a plane tree $T \in {\mathcal{LP}}_n$ ($n\geq 2$) is defined as follows:
\begin{itemize}
\item If the vertex  $\mathbf{j}$ is a singleton leaf, then label $\mathbf{j}$ by $a$;

\item If the vertex  $\mathbf{j}$ is an elder leaf, then label   $\mathbf{j}$ by $b$;

\item If the vertex  $\mathbf{j}$ is a young leaf, then label  $\mathbf{j}$ by $c$;
\item If the vertex  $\mathbf{j}$ is a young interior vertex, then label  $\mathbf{j}$ by $d$;
\item If the edge $(\mathbf{i},\mathbf{j})$ is a young edge, then label it by $t$.
\end{itemize}
 This  labeling scheme of   labeled plane trees  is called
the $(a,b,c,d,t)$-labeling.

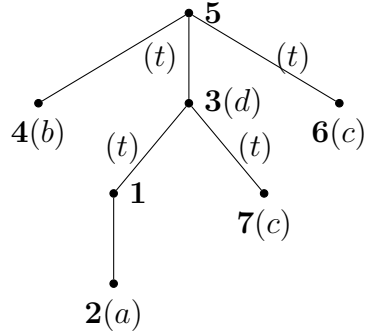
\begin{figure}[H]
    \centering
        \begin{tikzpicture}
[vertex/.style={shape=circle, draw, inner sep=1pt, fill=black},
sibling distance=2cm,level distance=1.2cm,
leaf/.style={label={[name=#1]below:$#1$}},auto]
\node[vertex,label=0:{$\mathbf{5}$}]{}[grow=down]
  child {node [vertex, label=-90:{$\mathbf{4}(b)$}]{}}
   child {node [vertex, label=0:{$\mathbf{3}(d)$}]{}
   child {node [vertex, label=0:{$\mathbf{1}$}]{}
   child {node [vertex, label=270:{$\mathbf{2}(a)$}]{}}edge from parent node[left]{$(t)$}}
   child {node [vertex, label=270:{$\mathbf{7}(c)$}]{}edge from parent node[right]{$(t)$}}edge from parent node[left]{$(t)$}}
    child {node [vertex, label=270:{$\mathbf{6}(c)$}]{}edge from parent node[right]{$(t)$}};
\end{tikzpicture}
     \caption{A plane tree in $\mathcal{LP}_6$ with the $(a,b,c,d,t)$-labeling.}
     \label{SOY-2}
     \end{figure}

For example, Fig. \ref{SOY-2} gives the grammatical labeling of  a  plane tree $T$ in ${\mathcal{LP}_{6}}$.  The labels of  $T$
 are shown in parentheses. 

The weight of $T$ is defined to be the product of all the labels, that is,
\[{\rm wt}(T)=a^{{\rm sleaf}(T)} b^{{\rm eleaf}(T)}c^{{\rm yleaf}(T)} d^{{\rm yint}(T)}t^{{\rm yedge}(T)}.\]
In this case, the weight of $T$  is ${\rm wt}(T)=t^4abc^2d$.

From the definition of the $(a,b,c,d,t)$-labeling, we see that  polynomial $\widetilde{G}_{n}(a,b,c,d;t)$ given in \eqref{defi-G} can be interpreted as: 
\begin{equation}
(n+1)!\widetilde{G}_{n}(a,b,c,d;t)=\sum_{T\in \mathcal{LP}_{n}}{\rm wt}(T)
\end{equation}
and so the proof of Theorem \ref{thm-SOY-grammar} is equivalent to the proof of the following assertion: For $n\geq 2$, 
\begin{equation}\label{gram-induct}
D^{n}(d)=\sum_{T\in \mathcal{LP}_{n}}{\rm wt}(T),
\end{equation}
where $D$ is the formal derivative with respect to the grammar defined in \eqref{SOY-grammar}.

We are now in a position to prove Theorem \ref{thm-SOY-grammar}.

\begin{figure}[H]
    \centering
    \begin{tikzpicture}
\node[shape=circle, draw, inner sep=1pt,fill=black!70, label=0:{$\mathbf{i}(d)$}]{}[grow=down]
[sibling distance=15mm,level distance=12mm]
child {node [shape=circle, draw, inner sep=1pt, fill=black!70, label=0:{$\mathbf{j}$}]{}[level distance=12mm]
child {node [shape=circle, draw, inner sep=1pt, fill=black!70, label=0:{$\mathbf{k}(a)$}]{}}edge from parent node[right]{$(t)$}};
\end{tikzpicture}\qquad
\begin{tikzpicture}
\node[shape=circle, draw, inner sep=1pt,fill=black!70, label=0:{$\mathbf{i}$}]{}[grow=down]
[sibling distance=12mm,level distance=12mm]
child {node [shape=circle, draw, inner sep=1pt, fill=black!70, label=-90:{$\mathbf{j}(b)$}]{}}
child {node [shape=circle, draw, inner sep=1pt, fill=black!70, label=-90:{$\mathbf{k}(c)$}]{}edge from parent node[right]{$(t)$}};
\end{tikzpicture}
    \caption{Two plane trees in $\mathcal{LP}_{2}$ with the $(a,b,c,d,t)$-labeling.}
    \label{n=2}
\end{figure}
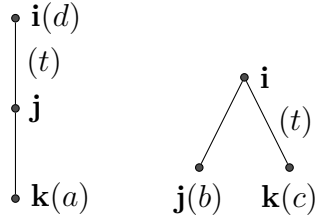

\noindent{\bf Proof of Theorem \ref{thm-SOY-grammar}.} We proceed by induction on $n$. When $n=2$, there are two plane trees with two edges and so there are twelve labeled plane trees in $\mathcal{LP}_{2}$. Fig. \ref{n=2} gives the grammatical labeling of  labeled plane trees in $\mathcal{LP}_{2}$ with  $\mathbf{i}$, $ \mathbf{j}$, $\mathbf{k}$ representing their vertices.

It is easy to check that 
\[\sum_{T\in \mathcal{LP}_{2}}{\rm wt}(T)=6tad+6tbc=D^2(d).\]
Hence this assertion is valid for $n=2$. Assume that it   holds for  $n-1$, to show that this assertion is also valid for $n$, by \eqref{gram-induct},  it suffices to show that 
\begin{equation}\label{gram-inductaa}
D\left(\sum_{T\in \mathcal{LP}_{n-1}}{\rm wt}(T)\right)=\sum_{T\in \mathcal{LP}_{n}}{\rm wt}(T).
\end{equation}
To this end, we will construct a  labeled plane trees $\widetilde{T}$ by inserting the vertex $\mathbf{n+1}$ into a  labeled plane tree   $T$ so that the change of weights between these two plane trees adheres to the substitution rules given in the grammar \eqref{SOY-grammar}. Given a plane tree $T \in \mathcal{LP}_{n-1}$, and let $\widetilde{T}$ be the  plane tree with $n+1$ vertices obtained by   performing a specific operation at position   $\diamond$ of $T$. We consider the following five cases:

\begin{itemize}
 \item[(1)] The position $\diamond$ is a singleton  leaf $\mathbf{j}$ of $T$ labeled by $a$. Suppose that  $\mathbf{i}$ is the parent of $\mathbf{j}$, see Fig. \ref{(1)}. We have two ways to create the plane tree $\widetilde{T}$.
 
  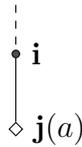
\begin{figure}[H]
    \centering
    \begin{tikzpicture}
\node{}[grow=down]
[sibling distance=15mm,level distance=8mm]
child {node [shape=circle, draw, inner sep=1pt, fill=black!70, label=0:{$\mathbf{i}$}]{}[level distance=10mm]
child {node [diamond,  draw, inner sep=1.3pt, fill=white!70, solid, label=0:{$\mathbf{j}(a)$}]{}edge from parent [solid]} edge from parent [dashed]};

\end{tikzpicture}
    \caption{The position
$\diamond$ is a singleton leaf of $T$.}
    \label{(1)}
\end{figure}

 \begin{itemize}
\item[(1-i)] Insert the vertex  $\mathbf{n+1}$ as   illustrated in Fig. \ref{(1-i)} to create the plane tree $\widetilde{T}$. It is easy to check that  the change of weights follows the substitution rule  $a\rightarrow 3tad$.

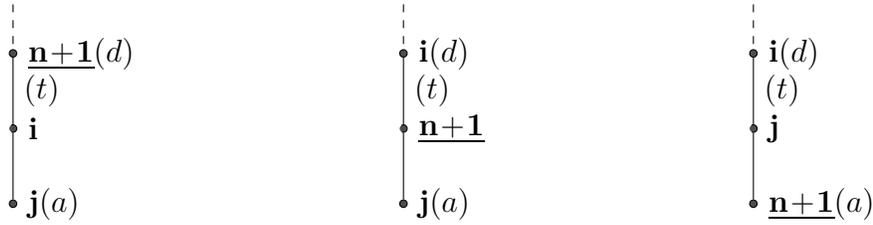
\begin{figure}[H]
    \centering
    \begin{subfigure}[t]{0.3\textwidth}
    \centering 
    \begin{tikzpicture}

\node{}[grow=down]
[sibling distance=15mm,level distance=8mm]
  child {node [shape=circle, draw,  inner sep=1pt, fill=black!70, label=0:{$\uline{\mathbf{n\!+\!1}}(d)$}]{}[level distance=10mm]
   child {node [shape=circle, draw, inner sep=1pt, fill=black!70, label=0:{$\mathbf{i}$}]{}[level distance = 10mm]
     child {node [shape=circle, draw, inner sep=1pt, fill=black!70, label=0:{$\mathbf{j}(a)$}]{} } 
edge from parent [solid] node[right]{$(t)$}}
edge from parent [dashed]
};
\end{tikzpicture}
       
        \caption{ $\mathbf{n+1}$ is a young interior vertex of $\widetilde{T}$ that has only one child, and its child also has only one child.}
        \label{(1-i-1)}
    \end{subfigure}\quad
    \begin{subfigure}[t]{0.3\textwidth}
    \centering 
        \begin{tikzpicture}

\node{}[grow=down]
[sibling distance=15mm,level distance=8mm]
  child {node [shape=circle, draw, inner sep=1pt, fill=black!70, label=0:{$\mathbf{i}(d)$}]{}[level distance=10mm]
   child {node [shape=circle, draw, inner sep=1pt, fill=black!70, label=0:{$\uline{\mathbf{n\!+\!1}}$}]{}[level distance = 10mm]
     child {node [shape=circle, draw, inner sep=1pt, fill=black!70, label=0:{$\mathbf{j}(a)$}]{} } 
edge from parent [solid] node[right]{$(t)$}}
edge from parent [dashed]
};
\end{tikzpicture}
        \caption{$\mathbf{n+1}$ is the parent of a singleton leaf in $\widetilde{T}$ and $\mathbf{n+1}$ is the only child of its parent.}
        \label{(1-i-2)}
    \end{subfigure}\quad
    \begin{subfigure}[t]{0.3\textwidth}
    \centering 
         \begin{tikzpicture}

\node{}[grow=down]
[sibling distance=15mm,level distance=8mm]
  child {node [shape=circle, draw, inner sep=1pt, fill=black!70, label=0:{$\mathbf{i}(d)$}]{}[level distance=10mm]
   child {node [shape=circle, draw, inner sep=1pt, fill=black!70, label=0:{$\mathbf{j}$}]{}[level distance = 10mm]
     child {node [shape=circle, draw, inner sep=1pt, fill=black!70, label=0:{$\uline{\mathbf{n\!+\!1}}(a)$}]{} } 
edge from parent [solid] node[right]{$(t)$}}
edge from parent [dashed]
};
\end{tikzpicture}
        \caption{$\mathbf{n+1}$ is a singleton leaf of  $\widetilde{T}$ and its parent has no  siblings.}
        \label{(1-i-3)}
    \end{subfigure}
    \caption{$a\rightarrow 3tad$.}
    \label{(1-i)}
\end{figure}

\item[(1-ii)] Insert the vertex  $\mathbf{n+1}$ according to the illustration in  Fig. \ref{(1-ii)} to form the  plane tree $\widetilde{T}$.  
A simple check reveals that the weights change according to the substitution rule $a\rightarrow 3tbc$.

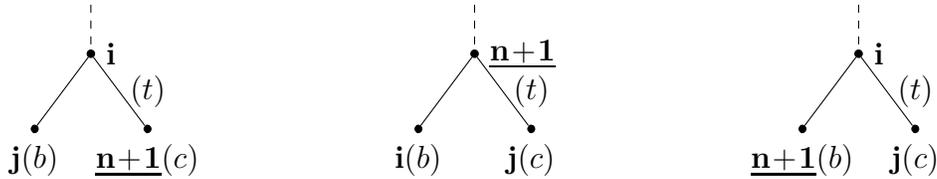
\begin{figure}[H]
    \centering
    \begin{subfigure}[t]{0.3\textwidth}
    \centering 
        \begin{tikzpicture}
[vertex/.style={shape=circle, draw, inner sep=1pt, fill=black},
sibling distance=1.5cm,level distance=0.8cm,
leaf/.style={label={[name=#1]below:$#1$}},auto]
\node{}[grow=down]
  child {node [vertex, label=0:{$\mathbf{i}$}]{}[level distance=10mm]
   child {node [vertex, label=-90:{$\mathbf{j}(b)$}]{}[level distance = 10mm]edge from parent [solid]}
     child {node [vertex, label=-90:{$\uline{\mathbf{n\!+\!1}}(c)$}]{}  
edge from parent [solid] node[right]{$(t)$}}
edge from parent [dashed]
};
\end{tikzpicture}
        \caption{$\mathbf{n+1}$ is a young leaf of $\widetilde{T}$ with only one sibling and its slbling is an elder leaf.}
        \label{(1-ii-1)}
    \end{subfigure}\quad
    \begin{subfigure}[t]{0.3\textwidth}
    \centering 
        \begin{tikzpicture}
[vertex/.style={shape=circle, draw, inner sep=1pt, fill=black},
sibling distance=1.5cm,level distance=0.8cm,
leaf/.style={label={[name=#1]below:$#1$}},auto]
\node{}[grow=down]
  child {node [vertex, label=0:{$\uline{\mathbf{n\!+\!1}}$}]{}[level distance=10mm]
   child {node [vertex, label=-90:{$\mathbf{i}(b)$}]{}[level distance = 10mm]edge from parent [solid]}
     child {node [vertex, label=-90:{$\mathbf{j}(c)$}]{}  
edge from parent [solid] node[right]{$(t)$}}edge from parent [dashed] };
\end{tikzpicture}
        \caption{$\mathbf{n+1}$ is the parent of an elder leaf  of $\widetilde{T}$, which 
 has only  two children, and the other child is a young leaf.}
        \label{(1-ii-2)}
    \end{subfigure}\quad
    \begin{subfigure}[t]{0.3\textwidth}
    \centering 
        \begin{tikzpicture}
[vertex/.style={shape=circle, draw, inner sep=1pt, fill=black},
sibling distance=1.5cm,level distance=0.8cm,
leaf/.style={label={[name=#1]below:$#1$}},auto]
\node{}[grow=down]
  child {node [vertex, label=0:{$\mathbf{i}$}]{}[level distance=10mm]
   child {node [vertex, label=-90:{$\uline{\mathbf{n\!+\!1}}(b)$}]{}[level distance = 10mm]edge from parent [solid]}
     child {node [vertex, label=-90:{$\mathbf{j}(c)$}]{}  
edge from parent [solid] node[right]{$(t)$}}edge from parent [dashed]};
\end{tikzpicture}
        \caption{$\mathbf{n+1}$ is an elder leaf of $\widetilde{T}$ with only one sibling and its slbling is a young leaf.}
        \label{(1-ii-3)}
    \end{subfigure}
    \caption{$a\rightarrow 3tbc$.}
    \label{(1-ii)}
\end{figure}
 
\end{itemize}

\item[(2)]  The position $\diamond$ is  an elder leaf $\mathbf{j}$ of $T$ labeled by $b$, see  Fig. \ref{(2)}. It follows the same actions as in Case(1) (see Fig. \ref{(2-i)} and Fig. \ref{(2-ii)}). It is straightforward to check that the resulting change in weights adheres to the substitution rule  $b\rightarrow 3t(ad+bc)$.

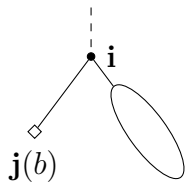
\begin{figure}[H]
    \centering
        \begin{tikzpicture}
[vertex/.style={shape=circle, draw, inner sep=1pt, fill=black},
subtree/.style={shape=ellipse, draw,minimum width=1.5cm, minimum height=.5cm},
sibling distance=1.5cm,level distance=0.8cm,
leaf/.style={label={[name=#1]below:$#1$}},auto]
\node{}[grow=down]
  child {node [vertex, label=0:{$\mathbf{i}$}]{}[level distance=10mm]
   child {node [diamond, draw, solid, inner sep=1.3pt, label=270:{$\mathbf{j}(b)$}]{} edge from parent [solid]
}
     child {node [subtree, solid, rotate=126]{}  edge from parent [solid]}
edge from parent [dashed]};
\end{tikzpicture}
        \caption{The position $\diamond$ is an elder leaf of $T$.}
        \label{(2)}
     \end{figure}

 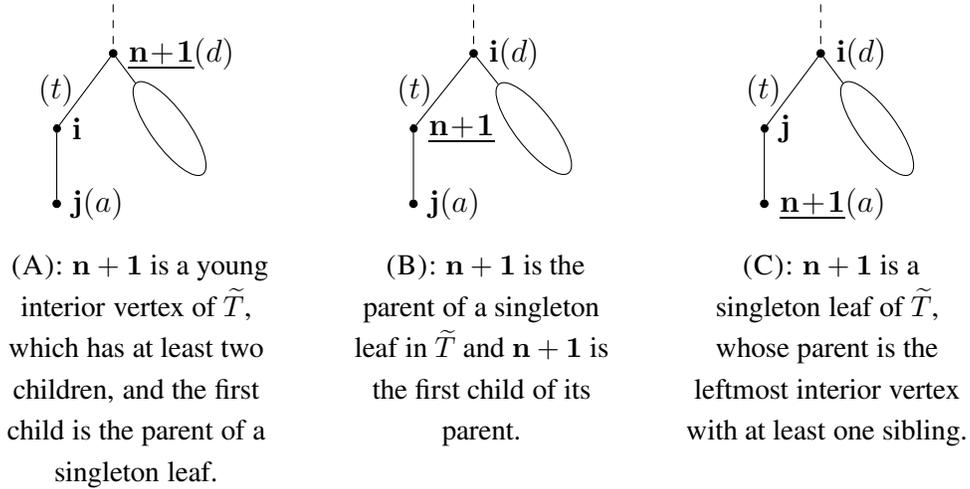
\begin{figure}[H]
    \centering
        \begin{subfigure}[t]{0.25\textwidth}
        \centering 
        \begin{tikzpicture}
[vertex/.style={shape=circle, draw, inner sep=1pt, fill=black},
subtree/.style={shape=ellipse, draw,minimum width=1.5cm, minimum height=.5cm},
sibling distance=1.5cm,level distance=0.8cm,
leaf/.style={label={[name=#1]below:$#1$}},auto]
\node{}[grow=down]
  child {node [vertex, label=0:{$\uline{\mathbf{n\!+\!1}}(d)$}]{}[level distance=10mm]
   child {node [vertex, label=0:{$\mathbf{i}$}]{}[level distance = 10mm]
   child {node [vertex, label=0:{$\mathbf{j}(a)$}]{}  }
edge from parent [solid] node[left]{$(t)$}}
     child {node [subtree,solid, rotate=126]{}  edge from parent [solid]}edge from parent [dashed]
};
\end{tikzpicture}
        \caption{ $\mathbf{n+1}$ is a young interior vertex of $\widetilde{T}$, which has at least two children, and the first child is the parent of a singleton leaf.}
        \label{(2-i-1)}
    \end{subfigure}\qquad
    \begin{subfigure}[t]{0.25\textwidth}
    \centering 
        \begin{tikzpicture}
[vertex/.style={shape=circle, draw, inner sep=1pt, fill=black},
subtree/.style={shape=ellipse, draw,minimum width=1.5cm, minimum height=.5cm},
sibling distance=1.6cm,level distance=0.8cm,
leaf/.style={label={[name=#1]below:$#1$}},auto]
\node{}[grow=down]
  child {node [vertex, label=0:{$\mathbf{i}(d)$}]{}[level distance=10mm]
   child {node [vertex, label=0:{$\uline{\mathbf{n\!+\!1}}$}]{}[level distance = 10mm]
   child {node [vertex, label=0:{$\mathbf{j}(a)$}]{}  }
edge from parent [solid] node[left]{$(t)$}}
     child {node [subtree,solid, rotate=128]{}  edge from parent [solid]}edge from parent [dashed]
};
\end{tikzpicture}
        \caption{$\mathbf{n+1}$ is the parent of a singleton leaf in $\widetilde{T}$  and $\mathbf{n+1}$ is the first child  of its parent. }  
        \label{(2-i-2)}
    \end{subfigure}\qquad
     \begin{subfigure}[t]{0.25\textwidth}
     \centering 
        \begin{tikzpicture}
[vertex/.style={shape=circle, draw, inner sep=1pt, fill=black},
subtree/.style={shape=ellipse, draw,minimum width=1.5cm, minimum height=.5cm},
sibling distance=1.5cm,level distance=0.8cm,
leaf/.style={label={[name=#1]below:$#1$}},auto]
\node{}[grow=down]
  child {node [vertex, label=0:{$\mathbf{i}(d)$}]{}[level distance=10mm]
   child {node [vertex, label=0:{$\mathbf{j}$}]{}[level distance = 10mm]
   child {node [vertex, label=0:{$\uline{\mathbf{n\!+\!1}}(a)$}]{}  }
edge from parent [solid] node[left]{$(t)$}}
     child {node [subtree,solid, rotate=126]{}  edge from parent [solid]}edge from parent [dashed]
};
\end{tikzpicture}
        \caption{ $\mathbf{n+1}$ is a singleton leaf of $\widetilde{T}$, whose parent is the leftmost interior vertex with at least one sibling.}
        \label{(2-i-3)}
    \end{subfigure}

    \caption{$b\rightarrow 3tad$.}
    \label{(2-i)}
\end{figure}
 
\begin{figure}[H]
    \centering
    \begin{subfigure}[t]{0.3\textwidth}
        \begin{tikzpicture}
[vertex/.style={shape=circle, draw, inner sep=1pt, fill=black},
subtree/.style={shape=ellipse, draw,minimum width=1.5cm, minimum height=.5cm},
sibling distance=1.5cm,level distance=0.8cm,
leaf/.style={label={[name=#1]below:$#1$}},auto]
\node{}[grow=down]
  child {node [vertex, label=0:{$\mathbf{i}$}]{}[level distance=15mm]
   child {node [vertex, label=270:{$\mathbf{j}(b)$}]{}edge from parent [solid] }
   child {node [vertex, label=270:{$\uline{\mathbf{n\!+\!1}}(c)$}]{}  edge from parent [solid]node[right, xshift=-1mm]{$(t)$}}
   child {node [subtree,solid, rotate=128]{}  edge from parent [solid]}edge from parent [dashed]
};
\end{tikzpicture}
        \caption{$\mathbf{n+1}$ is a young leaf of  $\widetilde{T}$  with at least two siblings  and the sibling immediately preceding  $\mathbf{n+1}$ is an elder leaf.}
        \label{(2-ii-1)}
    \end{subfigure}\quad
   \begin{subfigure}[t]{0.3\textwidth}
   \centering 
        \begin{tikzpicture}
[vertex/.style={shape=circle, draw, inner sep=1pt, fill=black},
subtree/.style={shape=ellipse, draw,minimum width=1.5cm, minimum height=.5cm},
sibling distance=1.5cm,level distance=0.8cm,
leaf/.style={label={[name=#1]below:$#1$}},auto]
\node{}[grow=down]
  child {node [vertex, label=0:{$\uline{\mathbf{n\!+\!1}}$}]{}[level distance=15mm]
   child {node [vertex, label=270:{$\mathbf{i}(b)$}]{}edge from parent [solid] }
   child {node [vertex, label=270:{$\mathbf{j}(c)$}]{}  edge from parent [solid]node[right, xshift=-1mm]{$(t)$}}
   child {node [subtree,solid, rotate=128]{}  edge from parent [solid]}edge from parent [dashed]
};
\end{tikzpicture}
        \caption{$\mathbf{n+1}$ is  the parent of an elder leaf in $\widetilde{T}$, which has at least three children and the second child is a young leaf.}
        \label{(2-ii-2)}
    \end{subfigure}\quad
    \begin{subfigure}[t]{0.3\textwidth}
    \centering 
        \begin{tikzpicture}
[vertex/.style={shape=circle, draw, inner sep=1pt, fill=black},
subtree/.style={shape=ellipse, draw,minimum width=1.5cm, minimum height=.5cm},
sibling distance=1.5cm,level distance=0.8cm,
leaf/.style={label={[name=#1]below:$#1$}},auto]
\node{}[grow=down]
  child {node [vertex, label=0:{$\mathbf{i}$}]{}[level distance=15mm]
   child {node [vertex, label=270:{$\uline{\mathbf{n\!+\!1}}(b)$}]{}edge from parent [solid] }
   child {node [vertex, label=270:{$\mathbf{j}(c)$}]{}  edge from parent [solid]node[right, xshift=-1mm]{$(t)$}}
   child {node [subtree,solid, rotate=128]{}  edge from parent [solid]}edge from parent [dashed]
};
\end{tikzpicture}
        \caption{$\mathbf{n+1}$ is an elder leaf of $\widetilde{T}$ with at least two siblings, and the sibling immediately following $\mathbf{n+1}$ is a young leaf.}
        \label{(2-ii-3)}
    \end{subfigure}
     \caption{$b\rightarrow 3tbc$.}
     \label{(2-ii)}
\end{figure}
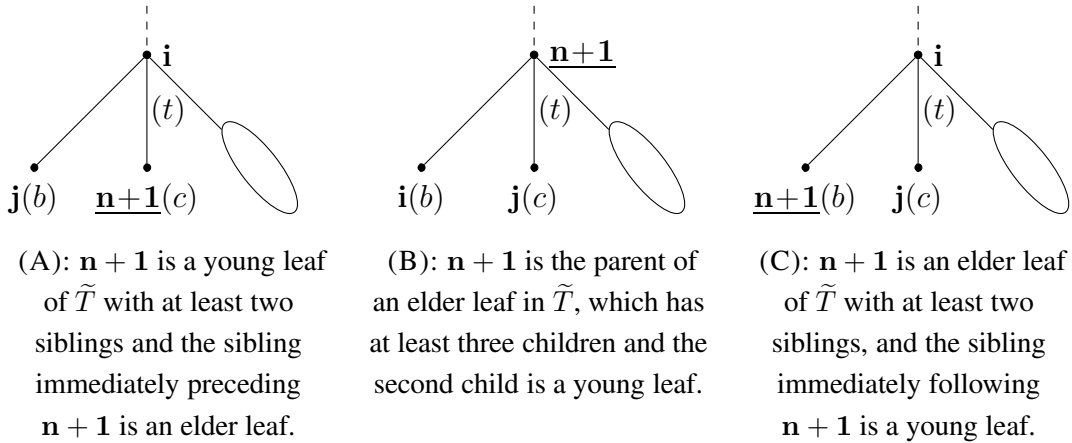

    \item[(3)] The position $\diamond$ is a young leaf  $\mathbf{j}$ of $T$ labeled by $c$,  as depicted in Fig. \ref{(3)}. We perform the steps shown in   Fig. \ref{(3-two)} to create the  plane tree $\widetilde{T}$. It is easy to check that   the change in weights follows the substitution rule $c\rightarrow 2a$.

\begin{figure}[H]
    \centering
     \begin{tikzpicture}
[vertex/.style={shape=circle, draw, inner sep=1pt, fill=black},
subtree/.style={shape=ellipse, draw,minimum width=1.5cm, minimum height=.5cm},
sibling distance=1.6cm,level distance=0.8cm,
leaf/.style={label={[name=#1]below:$#1$}},auto]
\node{}[grow=down]
  child {node [vertex, label=4:{$\mathbf{i}$}]{}[level distance=12mm]
   child {node [subtree, solid,rotate=-140]{} edge from parent [solid]}
   child {node [diamond, draw, solid, inner sep=1.3pt,label=270:{$\mathbf{j}(c)$}]{}edge from parent [solid]}
   child {node [subtree,solid, rotate=140]{} edge from parent [solid]}
edge from parent [dashed]};
\end{tikzpicture}
    \caption{The position $\diamond$ is a young leaf of $T$.}
    \label{(3)}
\end{figure}
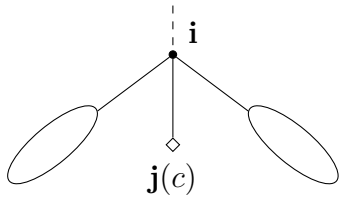
  
    \begin{figure}[H]
    \centering
    \begin{subfigure}[t]{0.4\textwidth}
        \begin{tikzpicture}
[vertex/.style={shape=circle, draw, inner sep=1pt, fill=black},
subtree/.style={shape=ellipse, draw,minimum width=1.5cm, minimum height=.5cm},
sibling distance=1.7cm,level distance=0.8cm,
leaf/.style={label={[name=#1]below:$#1$}},auto]
\node{}[grow=down]
  child {node [vertex, label=4:{$\mathbf{i}$}]{}[level distance=10mm]
   child {node [subtree, solid,rotate=-145]{} edge from parent [solid]}
   child {node [vertex, label=0:{$\mathbf{j}$}]{}
   child {node [vertex, label=0:{$\uline{\mathbf{n\!+\!1}}(a)$}]{}}
   edge from parent [solid] }
   child {node [subtree, solid,rotate=145]{} edge from parent [solid]}edge from parent [dashed]
};
\end{tikzpicture}
        \caption{$\mathbf{n+1}$ is a singleton leaf  of $\widetilde{T}$ whose parent has at least one elder sibling.}
        \label{(3-i)}
    \end{subfigure}\qquad
   \begin{subfigure}[t]{0.4\textwidth}
        \begin{tikzpicture}
[vertex/.style={shape=circle, draw, inner sep=1pt, fill=black},
subtree/.style={shape=ellipse, draw,minimum width=1.5cm, minimum height=.5cm},
sibling distance=1.7cm,level distance=0.8cm,
leaf/.style={label={[name=#1]below:$#1$}},auto]
\node{}[grow=down]
  child {node [vertex, label=4:{$\mathbf{i}$}]{}[level distance=10mm]
   child {node [subtree, solid,rotate=-145]{} edge from parent [solid]}
   child {node [vertex, label=0:{$\uline{\mathbf{n\!+\!1}}$}]{}
   child {node [vertex, label=0:{$\mathbf{j}(a)$}]{}}
   edge from parent [solid] }
   child {node [subtree, solid,rotate=145]{} edge from parent [solid]}edge from parent [dashed]
};
\end{tikzpicture}
\caption{$\mathbf{n+1}$ is the parent of a singleton leaf in $\widetilde{T}$  and  $\mathbf{n+1}$  is not the first child of its parent.}
        \label{(3-ii)}
    \end{subfigure}
     \caption{$c\rightarrow 2a$.}
     \label{(3-two)}
\end{figure}
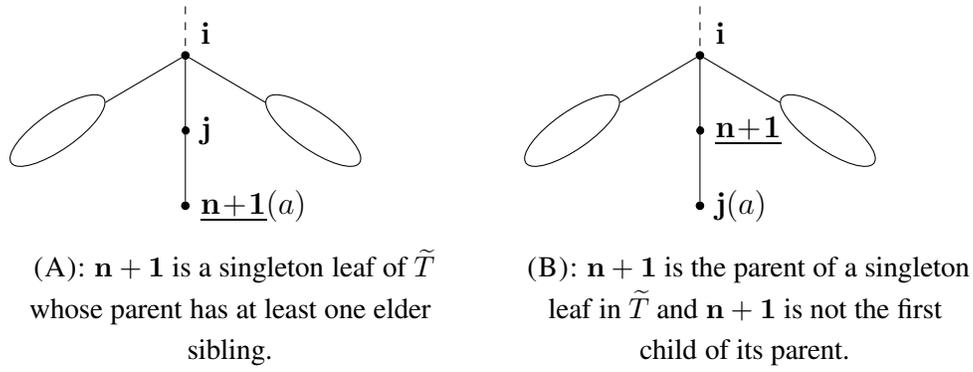

\item[(4)] The position $\diamond$ is depicted in Fig. \ref{(4)} as a young interior vertex $\mathbf{i}$ of $T$ labeled by $d$.  We execute the steps shown in Fig. \ref{(4-two)} to generate  the  plane tree $\widetilde{T}$. We see that the change in weights complies with  the substitution rule $d\rightarrow 2b$.

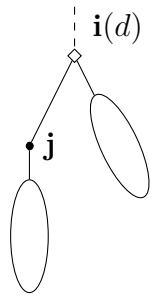
\begin{figure}[H]
    \centering
     \begin{tikzpicture}
[vertex/.style={shape=circle, draw, inner sep=1pt, fill=black},
subtree/.style={shape=ellipse, draw,minimum width=1.5cm, minimum height=.5cm},
sibling distance=1.2cm,level distance=0.8cm,
leaf/.style={label={[name=#1]below:$#1$}},auto]
\node{}[grow=down]
  child {node [diamond, draw, solid, inner sep=1.3pt, label=4:{$\mathbf{i}(d)$}]{}[level distance=12mm]
   child {node [vertex, label=0:{$\mathbf{j}$}]{}
   child {node [subtree, rotate=90]{}}
  edge from parent [solid]}
   child {node [subtree, solid,rotate=-65]{} edge from parent [solid]}
edge from parent [dashed]};
\end{tikzpicture}
        \caption{The position $\diamond$ is a young interior vertex of $T$.}
        \label{(4)}
\end{figure} 

\begin{figure}[H]
    \centering
    \begin{subfigure}[t]{0.3\textwidth}
        \begin{tikzpicture}
[vertex/.style={shape=circle, draw, inner sep=1pt, fill=black},
subtree/.style={shape=ellipse, draw,minimum width=1.5cm, minimum height=.5cm},
sibling distance=1.2cm,level distance=0.8cm,
leaf/.style={label={[name=#1]below:$#1$}},auto]
\node{}[grow=down]
  child {node [vertex, label=4:{$\mathbf{i}$}]{}[level distance=12mm]
  child {node [vertex, label=270:{$\uline{\mathbf{n\!+\!1}}(b)$}]{}edge from parent [solid]}
   child {node [vertex, label=0:{$\mathbf{j}$}]{}
   child {node [subtree, rotate=90]{}}
  edge from parent [solid]}
   child {node [subtree, solid,rotate=135]{} edge from parent [solid]}edge from parent [dashed]
};
\end{tikzpicture}
        \caption{$\mathbf{n+1}$ is an elder leaf of $\widetilde{T}$ with at least one sibling and the sibling immediately following   $\mathbf{n+1}$ is not a leaf.}
        \label{(4-i)}
    \end{subfigure}\qquad
   \begin{subfigure}[t]{0.3\textwidth}
        \begin{tikzpicture}
[vertex/.style={shape=circle, draw, inner sep=1pt, fill=black},
subtree/.style={shape=ellipse, draw,minimum width=1.5cm, minimum height=.5cm},
sibling distance=1.2cm,level distance=8mm,
leaf/.style={label={[name=#1]below:$#1$}},auto]
\node{}[grow=down]
  child {node [vertex, label=4:{$\uline{\mathbf{n\!+\!1}}$}]{}[level distance=12mm]
  child {node [vertex, label=270:{$\mathbf{i}(b)$}]{} edge from parent [solid]}
   child {node [vertex, label=0:{$\mathbf{j}$}]{}
   child {node [subtree, rotate=90]{}}
 edge from parent [solid]}
   child {node [subtree, solid,rotate=135]{} edge from parent [solid]}edge from parent [dashed]
};
\end{tikzpicture}
        \caption{$\mathbf{n+1}$  is the parent of  an elder leaf in $\widetilde{T}$, whose second child is not a leaf.}
        \label{(4-ii)}
    \end{subfigure}
     \caption{$d\rightarrow 2b$.}
     \label{(4-two)}
\end{figure}
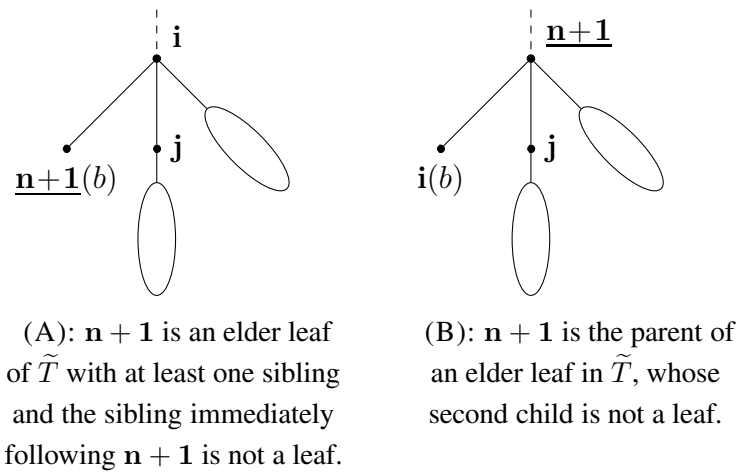

\item[(5)] The position $\diamond$ is  a young edge $(\mathbf{i},\mathbf{j})$ of $T$  labeled by $t$ as depicted in Fig. \ref{(5)}. We construct the plane tree $\widetilde{T}$ by following the steps shown in  Fig. \ref{(5-two)}. More precisely,  suppose $\mathbf{i}$ has $q$ children, which are   $\mathbf{k}_{1},\ldots,\mathbf{k}_{t-1},\mathbf{j},\mathbf{k}_{t+1},\ldots,\mathbf{k}_{q}$. We either add the edge $(\mathbf{i},\mathbf{n+1})$  immediately after $(\mathbf{i},\mathbf{j})$ or make   $\mathbf{n+1}$  the parent of $\mathbf{i}$.  In the latter case,  let $\mathbf{i},\mathbf{k}_{t+1},\ldots,\mathbf{k}_{q}$  be the successors of the vertex $\mathbf{n+1}$ (if exists) and let $\mathbf{k}_1,\ldots,\mathbf{k}_{t-1}$ and $\mathbf{j}$ be the children of $\mathbf{i}$ (if exists). It can be checked that the change in weights conforms to the substitution $t\rightarrow t^2(c+d)$.

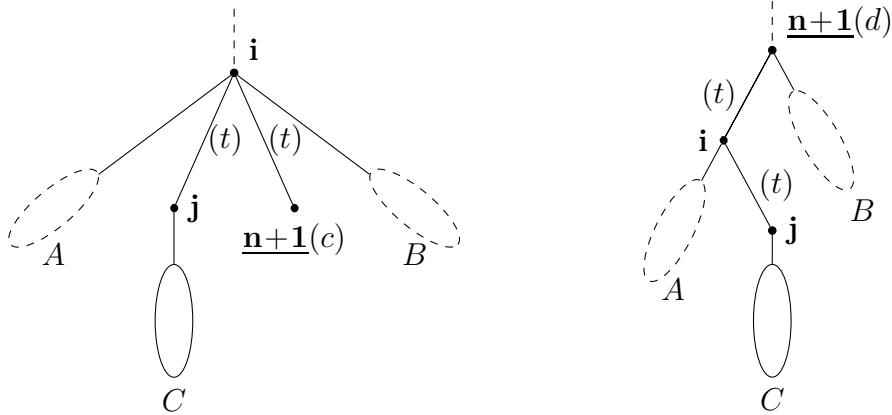
\begin{figure}[H]
    \centering
      \begin{tikzpicture}
[vertex/.style={shape=circle, draw, inner sep=1pt, fill=black},
subtree/.style={shape=ellipse, draw,minimum width=1.5cm, minimum height=.5cm},
every fit/.style={ellipse,draw,inner sep=-2pt},
sibling distance=1.5cm,level distance=1cm,
leaf/.style={label={[name=#1]below:$ $}}]

\node{}[grow=down]
  child {node [vertex, label=30:{$\mathbf{i}$}]{}[level distance=12mm]
  child {node [vertex, label=180:{$\mathbf{k}_1$}] (a's parent) {}  [level distance=12mm]
   child {node [subtree, rotate=90,leaf=a,label=-90:{$\cdots$}]{} [level distance=9mm] child {node [ rotate=70,leaf=c]{} edge from parent [solid]}edge from parent [solid]}edge from parent [solid]}
   child {node [vertex, label=180:{$\mathbf{k}_{t-1}$}] (b's parent) {}  [level distance=12mm]
   child {node [subtree, rotate=90,leaf=b]{} edge from parent [solid]}edge from parent [solid]}   
   child {node [vertex, label=0:{$\mathbf{j}$}]{} [level distance=16mm]
   child {node [subtree, rotate=90,label=180:{$C$},shift={(0.4,0)}]{}edge from parent [solid]} edge from parent[solid] node[xshift=2.4mm,yshift=-1.4mm]{$\diamond(t)$}}   
   child {node [vertex, label=0:{$\mathbf{k}_{t+1}$}] (d's parent) {}  [level distance=12mm]
   child {node [subtree, rotate=90,leaf=d,label=-90:{$\cdots$}]{}[level distance=9mm] child {node [ rotate=70,leaf=f]{}edge from parent [solid] }edge from parent [solid] }edge from parent [solid]}
   child {node [vertex, label=0:{$\mathbf{k}_q$}] (e's parent) {}  [level distance=12mm]
   child {node [subtree, rotate=90,leaf=e]{} edge from parent [solid]}edge from parent [solid]}
   edge from parent [dashed]};
   \node[draw,dashed,xshift=-2mm,fit=(a's parent)(a) (b) (c)  ,label=below:$A$] {}; 
     \node [dashed,xshift=-2mm,fit=(e) (f) (d's parent),label=below:$B$] {};
\end{tikzpicture}
        \caption{The position $\diamond$ is a young edge of $T$.}
        \label{(5)}
\end{figure}
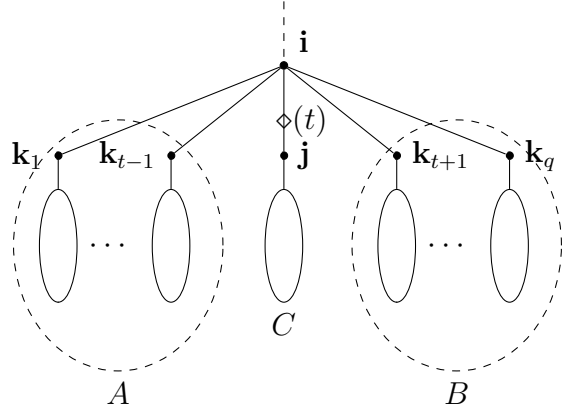

\begin{figure}[H]
    \centering
    \begin{subfigure}[t]{0.42\textwidth}
    \centering
        \begin{tikzpicture}
[vertex/.style={shape=circle, draw, inner sep=1pt, fill=black},
subtree/.style={shape=ellipse, draw,minimum width=1.5cm, minimum height=.5cm},
every fit/.style={ellipse,draw,inner sep=-2pt},
sibling distance=1.6cm,level distance=1cm,
leaf/.style={label={[name=#1]below:$ $}},auto]

\node{}[grow=down]
  child {node [vertex, label=30:{$\mathbf{i}$}]{}[level distance=18mm]
  child {node [subtree,rotate=40, label=230:{$A$}]{}edge from parent [solid]}
  child{node[vertex,label=0:{$\mathbf{j}$}] {}
  [level distance=15mm]
  child {node [subtree, rotate=90,label=180:{$C$}]{}edge from parent [solid]} edge from parent [solid]node[right,xshift=-1mm]{$(t)$}} 
   child {node [vertex,label=270:{$\uline{\mathbf{n\!+\!1}}(c)$} ]{}  edge from parent[solid]node[right,xshift=-1mm]{$(t)$}} 
   child {node [subtree, rotate=-40,label=310:{$B$}]{}edge from parent [solid]} 
   edge from parent [dashed]};
\end{tikzpicture}
        \caption{$\mathbf{n+1}$ is a young leaf of $\widetilde{T}$ with at least one sibling and the sibling immediately preceding $\mathbf{n+1}$ is not an elder leaf.}
        \label{(5-i)}
    \end{subfigure}\qquad
   \begin{subfigure}[t]{0.42\textwidth}
   \centering
        \begin{tikzpicture}
[vertex/.style={shape=circle, draw, inner sep=1pt, fill=black},
subtree/.style={shape=ellipse, draw,minimum width=1.5cm, minimum height=.5cm},
sibling distance=1.3cm,level distance=0.8cm,
leaf/.style={label={[name=#1]below:$#1$}},auto]
\node{}[grow=down]
  child {node [vertex, label=10:{$\uline{\mathbf{n\!+\!1}}(d)$}]{}[level distance=12mm]
  child {node [vertex, label=180:{$\mathbf{i}$}]{}
   child {node [subtree,  draw, dashed , rotate=63,label=205:{$A$}]{} edge from parent [solid]}
   child {node [vertex, label=0:{$\mathbf{j}$}]{}
   child {node [subtree,rotate=90,label=180:{$C$}]{}}edge from parent node[right]{$(t)$}}edge from parent node[left]{$(t)$}edge from parent [solid]}
   child {node [subtree,  draw, dashed ,rotate=120,label=170:{$B$}]{} edge from parent [solid]}edge from parent [dashed]
};
\end{tikzpicture}
        \caption{$\mathbf{n+1}$ is a young  interior vertex  of $\widetilde{T}$, whose first child is not the parent of a singleton leaf.}
        \label{(5-ii)}
    \end{subfigure}
     \caption{$t\rightarrow t^2(c+d)$.}
     \label{(5-two)}
\end{figure}

\end{itemize}

From the above construction, we find that 
\begin{itemize}
     \item  $\mathbf{n+1}$ is a singleton leaf of $\widetilde{T}$ in Fig. \ref{(1-i-3)}, Fig. \ref{(2-i-3)} and Fig. \ref{(3-i)} .
    \item   $\mathbf{n+1}$ is an elder leaf of $\widetilde{T}$ in Fig. \ref{(1-ii-3)}, Fig. \ref{(2-ii-3)} and Fig. \ref{(4-i)}.
    \item  $\mathbf{n+1}$ is a young leaf  of $\widetilde{T}$ in  Fig. \ref{(1-ii-1)}, Fig. \ref{(2-ii-1)} and Fig. \ref{(5-i)}.
     \item  $\mathbf{n+1}$ is the parent of a singleton leaf of $\widetilde{T}$ in Fig. \ref{(1-i-2)}, Fig. \ref{(2-i-2)} and Fig. \ref{(3-ii)}.
      \item  $\mathbf{n+1}$ is the parent of an elder leaf of $\widetilde{T}$ in Fig. \ref{(1-ii-2)}, Fig. \ref{(2-ii-2)} and Fig. \ref{(4-ii)}.
       \item  $\mathbf{n+1}$ is a young interior vertex of $\widetilde{T}$ in Fig. \ref{(1-i-1)}, Fig. \ref{(2-i-1)} and Fig. \ref{(5-ii)}.
\end{itemize}
Given these observations,  it is clear that this construction is  reversible according to the position of the vertex $\mathbf{n+1}$ in $\widetilde{T}$. 
Therefore, we could generate 
all labeled plane trees in $\mathcal{LP}_{n}$ from labeled plane trees in $\mathcal{LP}_{n-1}$. Moreover,  the change of weights between  plane trees in $\mathcal{LP}_{n}$  and  plane trees in $\mathcal{LP}_{n-1}$ in this construction is consistent with the substitution rule given in the grammar \eqref{SOY-grammar}.
Thus \eqref{gram-inductaa} is valid and consequently, the assertion is  valid for $n$ as well. This completes the proof.  \qed

\subsection{ The grammatical  derivation for Theorem  \ref{gf-gwilde} } For the formal derivative $D$ with respect to the grammar $G$ in \eqref{SOY-grammar}, we find that
$$
\begin{aligned}
& D(t^{-1})=-(c+d),\\[5pt]
& D^2(t^{-1})=-2(a+b),\\[5pt]
& D^3(t^{-1})=-2D(2b)=-2D^2(d).
\end{aligned}
$$
Consequently, invoking Theorem \ref{thm-SOY-grammar}, we derive that for $n\geq 3$, 
\[D^{n}(t^{-1})=-2n!\widetilde{G}_{n-1}(a,b,c,d,t).\]
Under the assumption $\widetilde{G}_{0}(a,b,c,d,t)=d$ and $\widetilde{G}_{1}(a,b,c,d,t)=b$, we deduce that  
\begin{align} \label{ttddeeqq}
{\rm Gen}(t^{-1};q)&=\sum_{n\geq 0} D^{n}(t^{-1})\frac{q^n}{n!} \nonumber \\[5pt]
&=t^{-1}-(c+d)q-(a+b)q^2+\sum_{n\geq 3} D^{n}(t^{-1})\frac{q^n}{n!} \nonumber \\[5pt]
&=t^{-1}-(c+d)q-(a+b)q^2-2\sum_{n\geq 3} n!\widetilde{G}_{n-1}(a,b,c,d,t)\frac{q^n}{n!}\nonumber \\[5pt]
&=t^{-1}-(c+d)q-(a+b)q^2-2q\sum_{n\geq 2}\widetilde{G}_{n}(a,b,c,d;t)q^n \nonumber \\[5pt]
&=t^{-1}-(c+d)q-(a+b)q^2-2q\sum_{n\geq 0}\widetilde{G}_{n}(a,b,c,d;t)q^n+2dq+2bq^2 \nonumber \\[5pt]
&=t^{-1}+(d-c)q+(b-a)q^2-2q\sum_{n\geq 0}\widetilde{G}_{n}(a,b,c,d;t)q^n.
\end{align}
In order to derive the generating function for $\widetilde{G}_{n}(a,b,c,d;t)$, it suffices to compute the generating function ${\rm Gen}(t^{-1};q)$. In view of the following relation 
\begin{equation}\label{aatttrr}
{\rm Gen}(t^{-2};q)={\rm Gen}^2(t^{-1};q),
\end{equation}
it is enough to compute ${\rm Gen}(t^{-2};q)$. 

Notice that 
$$
\begin{aligned}
& D(t^{-2})=-2t^{-1}(c+d),\\[5pt]
& D^2(t^{-2})=-4t^{-1}(a+b)+2(c+d)^2,\\[5pt]
& D^3(t^{-2})=12(a-b)(c-d),\\[5pt]
& D^4(t^{-2})=24(a-b)^2,\\[5pt]
& D^5(t^{-2})=0.
\end{aligned}
$$
Hence, we acquire 
\begin{align}\label{aattt}
{\rm Gen}(t^{-2};q) &= \sum_{n\geq 0}D^{n}(t^{-2})\frac{q^n}{n!}\notag\\[5pt]
&=t^{-2}-2t^{-1}(c+d)q+((c+d)^2-2t^{-1}(a+b))q^2 \notag\\[5pt]
&\quad+ 2(a-b)(c-d)q^3+(a-b)^2q^4.
\end{align}
By setting $q=0$ in \eqref{aatttrr}, we derive from \eqref{ttddeeqq} and \eqref{aattt} that  
\begin{align}\label{ttccdd}
&{\rm Gen}(t^{-1};q)=  \sqrt{{\rm Gen}(t^{-2};q)},
\end{align}
which completes the proof with the aid of  \eqref{ttddeeqq}, \eqref{aattt} and \eqref{ttccdd}.    \qed

\section{Proof of Theorem \ref{thm-gf-overlineM}  }

 To prove Theorem \ref{thm-gf-overlineM}  by using the grammatical approach, it is essential to  establish the grammar for the polynomial $M_n(u_1,u_2,u_3;v_1,v_2)$.  Additionally,  the concept of   the young edge introduced in Section 2 proves to be necessary. Recall that  ${\rm yedge}(T)$ denotes the number of young edges of a plane tree, i.e., the edges that are not connected to an old leaf. For $n\geq 1$, we define 
\begin{align*}
{M_n}(u_1,u_2,u_3;v_1,v_2;t)
&=\sum_{T\in\mathcal{T}_{n+1}} u_1^{{\rm sleaf}(T)}u_2^{\text {\rm etleaf}(T)}u_3^{\text {\rm entleaf}(T)}v_1^{\text{\rm yerleaf}(T)}v_2^{\text{\rm syleaf}(T)}t^{\text{\rm yedge}(T)}
\end{align*}
with the convention that 
${ M_0}(u_1,u_2,u_3;v_1,v_2;t)=u_3$. 

We first show that the following grammar 
\begin{align}\label{Refi-Motz}
M=\{& u_1\rightarrow 3tu_2v_2,~u_2\rightarrow 3tu_2v_1,~u_3\rightarrow 3tu_2v_2, \notag\\[4pt] &~ v_1\rightarrow 2(u_1+u_3),~v_2\rightarrow 4u_1u_2^{-1}u_3,~t\rightarrow t^2v_1\}
\end{align}
can be used to generate the polynomial ${M_n}(u_1,u_2,u_3;v_1,v_2;t)$. To wit, 

\begin{thm}\label{thm-Refi-Motz}
    Let $D$ be the formal derivative associated with the grammar \eqref{Refi-Motz}. For  $n\geq 1$, 
    \begin{equation}\label{thm-Refi-Motzee}
    D^n(2u_3)=(n+2)!{M_n}(u_1,u_2,u_3;v_1,v_2;t).
    \end{equation}
\end{thm}
The first few values of $D^n(2u_3)$ are given below: 
$$
\begin{aligned}
&D(2u_3) = 6tu_2v_2,\\[5pt]
&D^2(2u_3) = 24(t^2 u_2v_1v_2+ tu_1u_3),\\[5pt]
&D^3(2u_3)  = 120(t^3u_2v_1^2v_2+t^2u_1u_2v_2+t^2u_1u_3v_1  +t^2u_2u_3v_2),\\[5pt]
&D^4(2u_3)  = 720(t^4u_2v_1^3v_2 +t^3 u_2^2v_2^2+ 2t^3 u_2u_3v_1v_2 + 2t^3 u_1u_2v_1v_2 \\[5pt]
&\qquad\qquad + t^3u_1u_3v_1^2 + t^2u_1u_3^2 + t^2u_1^2u_3).
\end{aligned}
$$
Note that the grammar $M$ given in \eqref{gramar-Motzkin} for the Motzkin polynomial $M_n(u;v)$ can be derived  by  setting $v_1=v_2=v$ and $u_1=u_2=u_3=ut$ in Theorem \ref{thm-Refi-Motz}.  

Armed with Theorem \ref{thm-Refi-Motz}, we can provide a grammatical derivation of the generating function for ${M_n}(u_1,u_2,u_3;v_1,v_2;t)$. The special case $t=1$ then yields  the generating function for ${M_n}(u_1,u_2,u_3;v_1,v_2)$  as stated in Theorem  \ref{thm-gf-overlineM}. 

\begin{thm} \label{thm-gmotzkin} We have
    $$
\begin{aligned}
    &\quad\sum_{n\geq 0} {M_n}(u_1,u_2,u_3;v_1,v_2;t) q^{n}\\[5pt] 
    &= \frac{t^{-1}-v_1q+(u_3-u_1)q^2}{2q^2}\\[5pt]
    &\quad-\frac{\sqrt{t^{-2}-2t^{-1}v_1q+(v_1^2-2t^{-1}(u_1+u_3))q^2+ (2u_1v_1+2u_3v_1-4u_2v_2)q^3 +(u_1-u_3)^2q^4}}{2q^2}.
\end{aligned}
$$
\end{thm}

\subsection{A grammatical labeling of ${M_n}(u_1,u_2,u_3;v_1,v_2;t)$}  
\begin{figure}[H]
    \centering
        \begin{tikzpicture}
[vertex/.style={shape=circle, draw, inner sep=1pt, fill=black},
sibling distance=2cm,level distance=1.5cm]
\node[vertex,label=0:{$\mathbf{5}$}]{}[grow=down]
  child {node [vertex, label=-90:{$\mathbf{4}(u_3)$}]{}}
   child {node [vertex, label=0:{$\mathbf{3}$}]{}[sibling distance=1.5cm]
   child {node [vertex, label=-90:{$\mathbf{1}(u_2)$}]{}}
    child {node [vertex, label=-90:{$\mathbf{8}(v_2)$}]{}edge from parent node[right]{$(t)$}}
   child {node [vertex, label=0:{$\mathbf{7}$}]{}
    child {node [vertex, label=270:{$\mathbf{2}(u_1)$}]{}
   }edge from parent node[right]{$(t)$}}edge from parent node[right]{$(t)$}}
    child {node [vertex, label=270:{$\mathbf{6}(v_1)$}]{}edge from parent node[right]{$(t)$}};
\end{tikzpicture}
       \caption{A   tip-augmented plane tree  with the $(u_1,u_2,u_3;v_1,v_2;t)$-labeling.}
    \label{tipli-3}
     \end{figure}
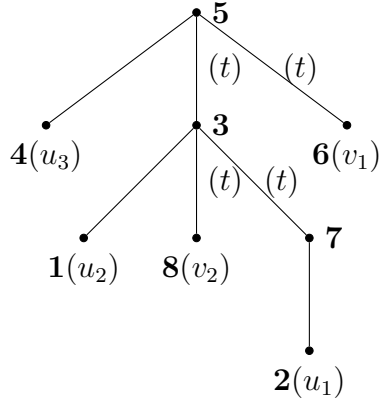

To prove Theorem \ref{thm-Refi-Motz}, we need to refine the grammatical labeling of   labeled  tip-augmented plane trees for the Motzkin polynomial $M_n(u;v)$ introduced in Section 2.  More precisely, 
 the refined grammatical labeling of a labeled tip-augmented plane tree $T\in{\mathcal{LT}_{n}}$ is defined as follows: 
\begin{itemize}
\item If the vertex $\mathbf{j}$ is a singleton leaf, then label $\mathbf{j}$ by $u_1$;

\item If the vertex  $\mathbf{j}$ is an elder twin leaf, then label  $\mathbf{j}$ by $u_2$;

\item If the vertex  $\mathbf{j}$ is  an elder non-twin leaf, then label $\mathbf{j}$ by  $u_3$;

\item If the vertex  $\mathbf{j}$ is a younger leaf, then label $\mathbf{j}$ by $v_1$;

\item If the vertex  $\mathbf{j}$ is a second leaf, then label $\mathbf{j}$ by $v_2$;

\item The young edge is labeled by $t$.
\end{itemize}
The weight of $T$ is defined to be the product of all the labels, that is, 
\[{\rm wt}(T)=u_1^{ {\rm sleaf}(T)} u_2^{\text {\rm etleaf}(T)}u_3^{\text {\rm entleaf}(T)}v_1^{\text{\rm yerleaf}(T)}v_2^{\text{\rm syleaf}(T)}t^{\text{\rm yedge}(T)}. \]
For example, Fig. \ref{tipli-3} shows  refined grammatical labeling of a   tip-augmented plane tree $T\in {\mathcal{LT}_7
}$ whose weight of $T$ is ${\rm wt}(T)=t^4u_1u_2u_3v_1v_2$.  We refer to this refined labeling of   tip-augmented plane trees as the $(u_1,u_2,u_3;v_1,v_2;t)$-labeling.

From the definition of the $(u_1,u_2,u_3;v_1,v_2;t)$-labeling, we see that the polynomial ${M_n}(u_1,u_2,u_3;v_1,v_2;t)$   can be interpreted as: 
\begin{equation}
(n+2)!{M}_{n}(u_1,u_2,u_3;v_1,v_2;t)=\sum_{T\in \mathcal{LT}_{n+1}}{\rm wt}(T)
\end{equation}

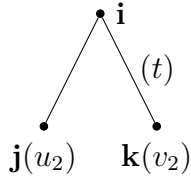
\begin{figure}[H]
    \centering
     \begin{tikzpicture}
[vertex/.style={shape=circle, draw, inner sep=1pt, fill=black},
sibling distance=1.5cm,level distance=1.5cm,]
\node[vertex, label=0:{$\mathbf{i}$}]{}[grow=down]
  child {node [vertex, label=270:{$\mathbf{j}(u_2)$}]{}}
  child {node [vertex, label=270:{$\mathbf{k}(v_2)$}]{}edge from parent node[right]{$(t)$}};
\end{tikzpicture}
    \caption{  A  tip-augmented plane tree $T \in \mathcal{LT}_2$  with the $(u_1,u_2,u_3;v_1,v_2;t)$-labeling.} \label{rfg-1}
\end{figure}

\noindent{\it Proof of Theorem \ref{thm-Refi-Motz}.} We proceed by induction on $n$. For $n=1$, there are six labeled tip-augmented plane trees $T$ in $\mathcal{LT}_{2}$. For every plane tree $T$ in  $\mathcal{LT}_{2}$ with $\mathbf{i},\mathbf{j},\mathbf{k}$ being its vertices, the $(u_1,u_2,u_3;v_1,v_2;t)$-labeling of $T$ is showed in Fig. \ref{rfg-1}, and  its weight is ${\rm wt}(T)=u_2v_2t$. As a result, it can be checked that \eqref{thm-Refi-Motzee} is valid for $n=1$. 
 Assume that this assertion holds for $n-1$. To show that it also holds for $n$,   it suffices to show that 
\begin{equation}\label{gram-inductaat}
D\left(\sum_{T\in \mathcal{LT}_{n}}{\rm wt}(T)\right)=\sum_{T\in \mathcal{LT}_{n+1}}{\rm wt}(T),
\end{equation}
where $D$ is the formal derivative with respect to the grammar defined in \eqref{Refi-Motz}. To accomplish this, we will construct   labeled tip-augmented plane tree in $\mathcal{LT}_{n+1}$ by adding the vertex $\mathbf{n+2}$ to labeled tip-augmented plane trees  in  $\mathcal{LT}_{n}$,  ensuring that the change in weights between these two plane trees adheres to the substitution rule given in the grammar \eqref{Refi-Motz}.

  Given a plane tree $T\in \mathcal{LT}_{n}$, and let $\widetilde{T}$ be the  plane tree with $n+2$ vertices obtained by  performing a specific operation   at the position $\diamond$ of $T$. We consider the following six cases: 
\begin{itemize}
 \item[(1).]  The position $\diamond$  is  a singleton leaf $\mathbf{j}$ of $T$  labeled by  $u_1$.  Suppose that $\mathbf{i}$ is the parent of $\mathbf{j}$, see Fig. \ref{(one)}.  Insert the vertex $\mathbf{n+2}$ as illustrated in Fig. \ref{(one-i)}. It is easy to check that the change of weights follows the substitution rule $u_1\rightarrow 3tu_2v_2$. 

 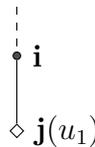
\begin{figure}[H]
    \centering
    \begin{tikzpicture}
\node{}[grow=down]
[sibling distance=15mm,level distance=8mm]
child {node [shape=circle, draw, inner sep=1pt, fill=black!70, label=0:{$\mathbf{i}$}]{}[level distance=10mm]
child {node [diamond,  draw, inner sep=1.3pt, fill=white!70, solid, label=0:{$\mathbf{j}(u_1)$}]{}edge from parent [solid]} edge from parent [dashed]};

\end{tikzpicture}
    \caption{The position
$\diamond$ is a singleton leaf of $T$.}
    \label{(one)}
\end{figure}

 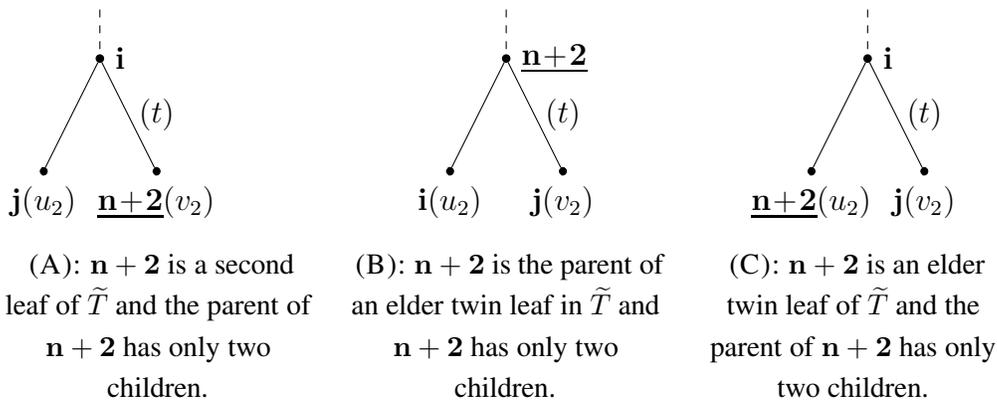
\begin{figure}[H]
    \centering
    \begin{subfigure}[t]{0.28\textwidth}
        \begin{tikzpicture}
[vertex/.style={shape=circle, draw, inner sep=1pt, fill=black},
subtree/.style={shape=ellipse, draw,minimum width=1.5cm, minimum height=.5cm},
sibling distance=1.5cm,level distance=0.8cm,
leaf/.style={label={[name=#1]below:$#1$}},auto]
\node{}[grow=down]
  child {node [vertex, label=0:{$\mathbf{i}$}]{}[level distance=15mm]
   child {node [vertex, label=270:{$\mathbf{j}(u_2)$}]{}edge from parent [solid] }
   child {node [vertex, label=270:{$\uline{\mathbf{n\!+\!2}}(v_2)$}]{}  edge from parent [solid] node[right]{$(t)$}}edge from parent [dashed]
};
\end{tikzpicture}
        \caption{ $\mathbf{n+2}$ is a second leaf of $\widetilde{T}$ and the parent of $\mathbf{n+2}$ has only two children. }
        \label{(one-1)}
    \end{subfigure}\quad
   \begin{subfigure}[t]{0.28\textwidth}
   \centering 
       \begin{tikzpicture}
[vertex/.style={shape=circle, draw, inner sep=1pt, fill=black},
subtree/.style={shape=ellipse, draw,minimum width=1.5cm, minimum height=.5cm},
sibling distance=1.5cm,level distance=0.8cm,
leaf/.style={label={[name=#1]below:$#1$}},auto]
\node{}[grow=down]
  child {node [vertex, label=0:{$\uline{\mathbf{n\!+\!2}}$}]{}[level distance=15mm]
   child {node [vertex, label=270:{$\mathbf{i}(u_2)$}]{}edge from parent [solid] }
   child {node [vertex, label=270:{$\mathbf{j}(v_2)$}]{}  edge from parent [solid] node[right]{$(t)$}}edge from parent [dashed]
};
\end{tikzpicture}
        \caption{$\mathbf{n+2}$ is  the parent of an  elder twin leaf in $\widetilde{T}$ and $\mathbf{n+2}$ has only two children. }
        \label{(one-2)}
    \end{subfigure}\quad
    \begin{subfigure}[t]{0.28\textwidth}
    \centering 
         \begin{tikzpicture}
[vertex/.style={shape=circle, draw, inner sep=1pt, fill=black},
subtree/.style={shape=ellipse, draw,minimum width=1.5cm, minimum height=.5cm},
sibling distance=1.5cm,level distance=0.8cm,
leaf/.style={label={[name=#1]below:$#1$}},auto]
\node{}[grow=down]
  child {node [vertex, label=0:{$\mathbf{i}$}]{}[level distance=15mm]
   child {node [vertex, label=270:{$\uline{\mathbf{n\!+\!2}}(u_2)$}]{}edge from parent [solid] }
   child {node [vertex, label=270:{$\mathbf{j}(v_2)$}]{}  edge from parent [solid] node[right]{$(t)$}}edge from parent [dashed]
};
\end{tikzpicture}
        \caption{$\mathbf{n+2}$ is an elder twin leaf of $\widetilde{T}$ and the parent of $\mathbf{n+2}$ has only two children.}
        \label{(one-3)}
    \end{subfigure}
     \caption{$u_1\rightarrow 3tu_2v_2$.}
     \label{(one-i)}
\end{figure}

    \item[(2).]  The position $\diamond$  is  an   elder twin leaf $\mathbf{j}$ of $T$ labeled by  $u_2$. Suppose that $\mathbf{k}$ is the twin sibling leaf of $\mathbf{j}$ and $\mathbf{i}$ is the parent of both $\mathbf{j}$ and $\mathbf{k}$, see Fig. \ref{(two)}. It follows the same actions as in Case (1), see Fig. \ref{(two-i)}. It is straightforward to check the resulting change in weights adheres to the substitution rule $u_2\rightarrow 3tu_2v_1$.

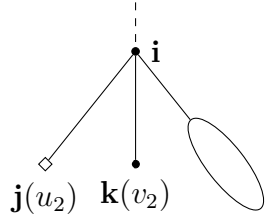
\begin{figure}[H]
    \centering
     \begin{tikzpicture}
[vertex/.style={shape=circle, draw, inner sep=1pt, fill=black},
subtree/.style={shape=ellipse, draw,minimum width=1.5cm, minimum height=.5cm},
sibling distance=1.2cm,level distance=0.8cm,
leaf/.style={label={[name=#1]below:$#1$}},auto]
\node{}[grow=down]
  child {node [vertex, label=0:{$\mathbf{i}$}]{}[level distance=15mm]
   child {node [diamond,draw, solid, inner sep=1.3pt, label=270:{$\mathbf{j}(u_2)$}]{}edge from parent [solid] }
   child {node [vertex, label=270:{$\mathbf{k}(v_2)$}]{}  edge from parent [solid]}
   child {node [subtree,solid, rotate=128]{}  edge from parent [solid]}edge from parent [dashed]
};
\end{tikzpicture}
 \caption{The position $\diamond$ is an elder twin leaf of $T$.}
    \label{(two)}
\end{figure}

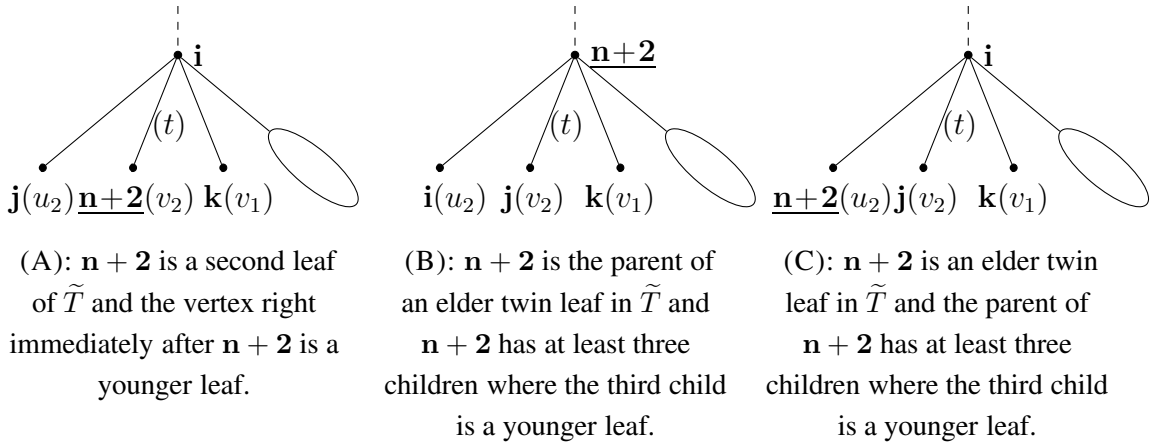
\begin{figure}[H]
    \centering
    \begin{subfigure}[t]{0.31\textwidth}
      \begin{tikzpicture}
[vertex/.style={shape=circle, draw, inner sep=1pt, fill=black},
subtree/.style={shape=ellipse, draw,minimum width=1.5cm, minimum height=.5cm},
sibling distance=1.2cm,level distance=0.8cm,
leaf/.style={label={[name=#1]below:$#1$}},auto]
\node{}[grow=down]
  child {node [vertex, label=0:{$\mathbf{i}$}]{}[level distance=15mm]
   child {node [vertex, label=270:{$\mathbf{j}(u_2)$}]{}edge from parent [solid] }
   child {node [vertex, label=270:{$~\uline{\mathbf{n\!+\! 2}}(v_2)$}]{}edge from parent [solid] node[right,xshift=-2mm,yshift=-2mm]{$(t)$} }
   child {node [vertex, label=270:{$\quad  \mathbf{k}(v_1)$}]{}  edge from parent [solid]}
   child {node [subtree,solid, rotate=140]{}  edge from parent [solid]}edge from parent [dashed]
};
\end{tikzpicture}
\caption{ $\mathbf{n+2}$ is a second leaf of $\widetilde{T}$ and the  vertex right immediately after  $\mathbf{n+2}$ is a younger leaf.}
\label{(two-1)}
    \end{subfigure}\quad
    \begin{subfigure}[t]{0.31\textwidth}
    \centering 
         \begin{tikzpicture}
[vertex/.style={shape=circle, draw, inner sep=1pt, fill=black},
subtree/.style={shape=ellipse, draw,minimum width=1.5cm, minimum height=.5cm},
sibling distance=1.2cm,level distance=0.8cm,
leaf/.style={label={[name=#1]below:$#1$}},auto]
\node{}[grow=down]
  child {node [vertex, label=0:{$\uline{\mathbf{n\!+\! 2}}$}]{}[level distance=15mm]
   child {node [vertex, label=270:{$\quad \mathbf{i}(u_2)$}]{}edge from parent [solid] }
   child {node [vertex, label=270:{$ ~\mathbf{j}(v_2)$}]{}edge from parent [solid] node[right,xshift=-2mm,yshift=-2mm]{$(t)$} }
   child {node [vertex, label=270:{$\mathbf{k}(v_1)$}]{}  edge from parent [solid]}
   child {node [subtree,solid, rotate=140]{}  edge from parent [solid]}edge from parent [dashed]
};
\end{tikzpicture}
        \caption{$\mathbf{n+2}$ is  the parent of an  elder twin leaf in $\widetilde{T}$ and $\mathbf{n+2}$ has at least three children where the third child is a younger leaf.}
        \label{(two-2)}
    \end{subfigure}\quad
    \begin{subfigure}[t]{0.31\textwidth}
   \centering 
         \begin{tikzpicture}
[vertex/.style={shape=circle, draw, inner sep=1pt, fill=black},
subtree/.style={shape=ellipse, draw,minimum width=1.5cm, minimum height=.5cm},
sibling distance=1.2cm,level distance=0.8cm,
leaf/.style={label={[name=#1]below:$#1$}},auto]
\node{}[grow=down]
  child {node [vertex, label=0:{$\mathbf{i}$}]{}[level distance=15mm]
   child {node [vertex, label=270:{$\uline{\mathbf{n\!+\! 2}}(u_2)$}]{}edge from parent [solid] }
   child {node [vertex, label=270:{$ ~\mathbf{j}(v_2)$}]{}edge from parent [solid] node[right,xshift=-2mm,yshift=-2mm]{$(t)$} }
   child {node [vertex, label=270:{$\mathbf{k}(v_1)$}]{}  edge from parent [solid]}
   child {node [subtree,solid, rotate=140]{}  edge from parent [solid]}edge from parent [dashed]
};
\end{tikzpicture}
        \caption{$\mathbf{n+2}$ is an elder twin leaf in $\widetilde{T}$ and the parent of $\mathbf{n+2}$ has at least three children where the third child is a younger leaf. }
        \label{(two-3)}
    \end{subfigure}
     \caption{$u_2\rightarrow 3tu_2v_1$.}
     \label{(two-i)}
\end{figure}
    
    \item[(3).]The position $\diamond$  is  an  elder non-twin leaf $\mathbf{j}$ of $T$ labeled by  $u_3$. Suppose that $\mathbf{k}$ is the sibling right immediately after $\mathbf{j}$ and    
    $\mathbf{i}$ is the parent of both $\mathbf{j}$ and $\mathbf{k}$, see Fig. \ref{(three)}. It takes the same performances as in Case(1), see Fig. \ref{(three-i)}. It is clear to find that the weights are satisfies with the substitute rule $u_3\rightarrow 3tu_2v_2$.

   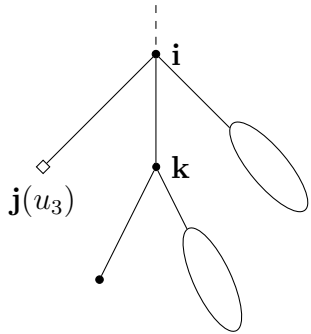

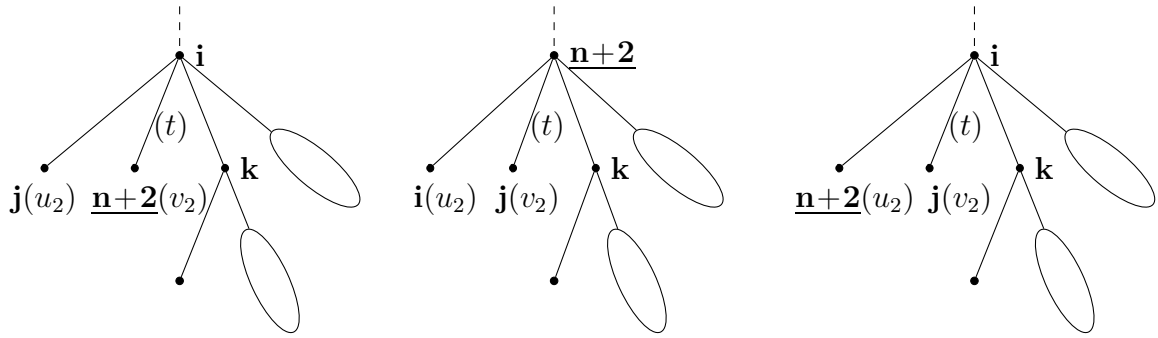
\begin{figure}[H]
    \centering
     \begin{tikzpicture}
[vertex/.style={shape=circle, draw, inner sep=1pt, fill=black},
subtree/.style={shape=ellipse, draw,minimum width=1.5cm, minimum height=.5cm},
sibling distance=1.5cm,level distance=0.8cm,
leaf/.style={label={[name=#1]below:$#1$}},auto]
\node{}[grow=down]
  child {node [vertex, label=0:{$\mathbf{i}$}]{}[level distance=15mm]
   child {node [diamond,draw, solid, inner sep=1.3pt, label=270:{$\mathbf{j}(u_3)$}]{}edge from parent [solid] }
   child {node [vertex, label=0:{$\mathbf{k}$}]{} 
   child {node[vertex] {} [level distance=10mm]
   edge from parent [solid]}
   child {node[subtree,rotate=115]{}
   edge from parent [solid]} edge from parent [solid]}
   child {node [subtree,solid, rotate=130]{}  edge from parent [solid]}edge from parent [dashed]
};
\end{tikzpicture}
 \caption{The position $\diamond$ is an elder non-twin leaf of $T$.}
    \label{(three)}
\end{figure}

\begin{figure}[H]
    \centering
    \begin{subfigure}[t]{0.3\textwidth}
     \begin{tikzpicture}
[vertex/.style={shape=circle, draw, inner sep=1pt, fill=black},
subtree/.style={shape=ellipse, draw,minimum width=1.5cm, minimum height=.5cm},
sibling distance=1.2cm,level distance=0.8cm,
leaf/.style={label={[name=#1]below:$#1$}},auto]
\node{}[grow=down]
  child {node [vertex, label=0:{$\mathbf{i}$}]{}[level distance=15mm]
   child {node [vertex, inner sep=1pt, label=270:{$\mathbf{j}(u_2)$}]{}edge from parent [solid] }
   child {node [vertex, inner sep=1pt, label=270:{$\quad \uline{\mathbf{n\! + \! 2}}(v_2)$}]{}edge from parent [solid]node[right,xshift=-2mm,yshift=-2mm]{$(t)$} }
   child {node [vertex, label=0:{$\mathbf{k}$}]{} 
   child {node[vertex] {} [level distance=10mm]
   edge from parent [solid]}
   child {node[subtree,rotate=115]{}
   edge from parent [solid]} edge from parent [solid]}
   child {node [subtree,solid, rotate=140]{}  edge from parent [solid]}edge from parent [dashed]
};
\end{tikzpicture}
\caption{$\mathbf{n+2}$ is a second leaf of $\widetilde{T}$ and the  vertex right immediately after  $\mathbf{n+2}$ is not a younger leaf.}
\label{(three-1)}
    \end{subfigure}\quad
     \begin{subfigure}[t]{0.31\textwidth}
    \centering 
       \begin{tikzpicture}
[vertex/.style={shape=circle, draw, inner sep=1pt, fill=black},
subtree/.style={shape=ellipse, draw,minimum width=1.5cm, minimum height=.5cm},
sibling distance=1.1cm,level distance=0.8cm,
leaf/.style={label={[name=#1]below:$#1$}},auto]
\node{}[grow=down]
  child {node [vertex, label=0:{$\uline{\mathbf{n\!+\!2}}$}]{}[level distance=15mm]
   child {node [vertex, inner sep=1pt, label=270:{$\quad \mathbf{i}(u_2)$}]{}edge from parent [solid] }
   child {node [vertex, inner sep=1pt, label=270:{$\quad \mathbf{j}(v_2)$}]{}edge from parent [solid]node[right,xshift=-2mm,yshift=-2mm]{$(t)$} }
   child {node [vertex, label=0:{$\mathbf{k}$}]{} 
   child {node[vertex] {} [level distance=10mm]
   edge from parent [solid]}
   child {node[subtree,rotate=115]{}
   edge from parent [solid]} edge from parent [solid]}
   child {node [subtree,solid, rotate=140]{}  edge from parent [solid]}edge from parent [dashed]
};
\end{tikzpicture}
        \caption{$\mathbf{n+2}$ is  the parent of an  elder twin leaf in $\widetilde{T}$ and $\mathbf{n+2}$ has at least three children where the third child is not a younger leaf.}
        \label{(three-2)}
    \end{subfigure}\quad
   \begin{subfigure}[t]{0.31\textwidth}
   \centering 
         \begin{tikzpicture}
[vertex/.style={shape=circle, draw, inner sep=1pt, fill=black},
subtree/.style={shape=ellipse, draw,minimum width=1.5cm, minimum height=.5cm},
sibling distance=1.2cm,level distance=0.8cm,
leaf/.style={label={[name=#1]below:$#1$}},auto]
\node{}[grow=down]
  child {node [vertex, label=0:{$\mathbf{i}$}]{}[level distance=15mm]
   child {node [vertex, inner sep=1pt, label=270:{$\quad \uline{\mathbf{n\!+\!2}}(u_2)$}]{}edge from parent [solid] }
   child {node [vertex, inner sep=1pt, label=270:{$\qquad \mathbf{j}(v_2)$}]{}edge from parent [solid]node[right,xshift=-2mm,yshift=-2mm]{$(t)$} }
   child {node [vertex, label=0:{$\mathbf{k}$}]{} 
   child {node[vertex] {} [level distance=10mm]
   edge from parent [solid]}
   child {node[subtree,rotate=115]{}
   edge from parent [solid]} edge from parent [solid]}
   child {node [subtree,solid, rotate=140]{}  edge from parent [solid]}edge from parent [dashed]
};
\end{tikzpicture}
        \caption{$\mathbf{n+2}$ is an elder twin leaf of $\widetilde{T}$ and the parent of $\mathbf{n+2}$ has at least three children where the third child is not a younger leaf.  }
        \label{(three-3)}
    \end{subfigure}\quad
   
     \caption{$u_3\rightarrow 3tu_2v_2$.}
     \label{(three-i)}
\end{figure}

\item[(4).] The position $\diamond$  is a younger leaf $\mathbf{j}$ of $T$ labeled by  $v_1$, as depicted in Fig. \ref{(four)}. It should be noted that the forest $A$,  formed  by the subtrees $\mathbf{k_{2}}, \ldots, \mathbf{k_{t-1}}$ can not be empty, while the forest $B$, formed by the subtrees $\mathbf{k_{t+1}}, \ldots, \mathbf{k_{q}}$ may be empty. 
\begin{itemize}
    \item[(4-i).] Insert the vertex $\mathbf{n+2}$ as illustrated in Fig. \ref{(four-i)} to form the tree $\widetilde{T}$. A simple check reveals the weights follows the substitution rule $v_1\rightarrow 2u_1$.
     \item[(4-ii).] Insert the vertex $\mathbf{n+2}$ as illustrated in Fig. \ref{(four-ii)} to construct the tree $\widetilde{T}$. A simple check reveals the weights follows the substitution rule $v_1\rightarrow 2u_3$.
\end{itemize}

\begin{figure}[H]
    \centering
      \begin{tikzpicture}
[vertex/.style={shape=circle, draw, inner sep=1pt, fill=black},
subtree/.style={shape=ellipse, draw,minimum width=1.5cm, minimum height=.5cm},
every fit/.style={ellipse,draw,inner sep=-2pt},
sibling distance=1.5cm,level distance=1cm,
leaf/.style={label={[name=#1]below:$ $}}]

\node{}[grow=down]
  child {node [vertex, label=30:{$\mathbf{i}$}]{}[level distance=12mm]
  child {node [vertex,label=180:{$\mathbf{k}_1$}]  {} edge from parent [solid]}
  child {node [vertex, label=180:{$\mathbf{k}_2$}] (a's parent) {}  [level distance=12mm]
   child {node [subtree, rotate=90,leaf=a,label=-90:{$\cdots$}]{} [level distance=9mm] child {node [ rotate=70,leaf=c]{} edge from parent [solid]}edge from parent [solid]}edge from parent [solid]}
   child {node [vertex, label=180:{$\mathbf{k}_{t-1}$}] (b's parent) {}  [level distance=12mm]
   child {node [subtree, rotate=90,leaf=b]{} edge from parent [solid]}edge from parent [solid]}   
   child {node [diamond, draw, solid, inner sep=1.3pt,label=0:{$\!\mathbf{j}(v_1)$}]{} edge from parent [solid]} 
   child {node [vertex, label=0:{$\mathbf{k}_{t+1}$}] (d's parent) {}  [level distance=12mm]
   child {node [subtree, rotate=90,leaf=d,label=-90:{$\cdots$}]{}[level distance=9mm] child {node [ rotate=70,leaf=f]{}edge from parent [solid] }edge from parent [solid] }edge from parent [solid]}
   child {node [vertex, label=0:{$\mathbf{k}_q$}] (e's parent) {}  [level distance=12mm]
   child {node [subtree, rotate=90,leaf=e]{} edge from parent [solid]}edge from parent [solid]}
   edge from parent [dashed]};
   \node[draw,dashed,xshift=-2mm,fit=(a's parent)(a) (b) (c)  ,label=below:$A$] {}; 
     \node [dashed,xshift=-2mm,fit=(e) (f) (d's parent),label=below:$B$] {};
\end{tikzpicture}
        \caption{The position $\diamond$ is a younger leaf of $T$.}
        \label{(four)}
\end{figure}
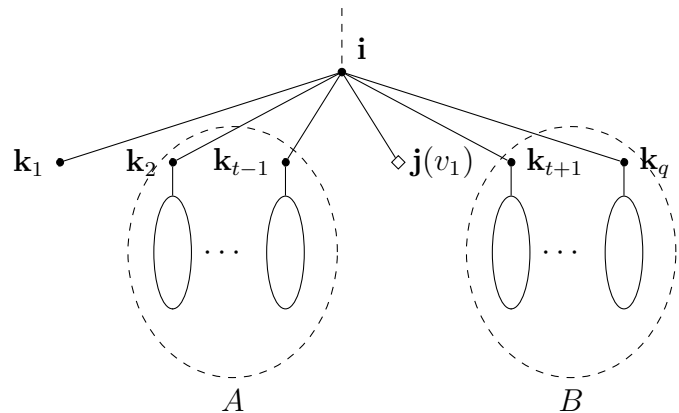


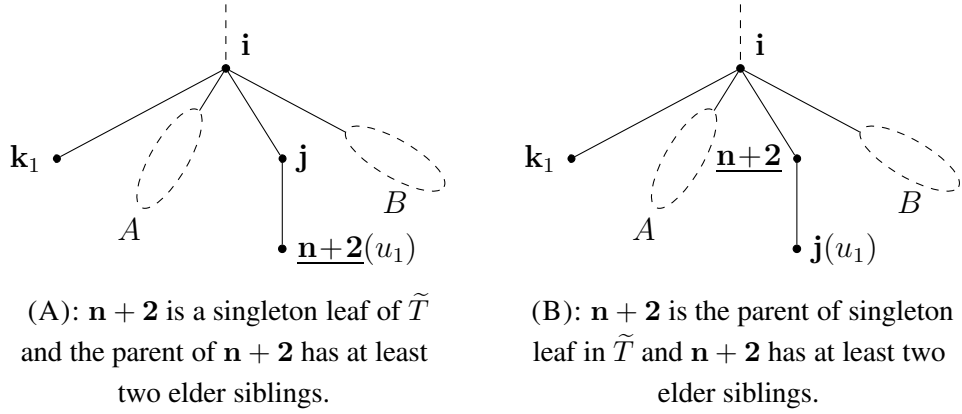
\begin{figure}[H]
    \centering
    \begin{subfigure}[t]{0.4\textwidth}
    \centering
        \begin{tikzpicture}
[vertex/.style={shape=circle, draw, inner sep=1pt, fill=black},
subtree/.style={shape=ellipse, draw,minimum width=1.5cm, minimum height=.5cm},
every fit/.style={ellipse,draw,inner sep=-2pt},
sibling distance=1.5cm,level distance=1cm,
leaf/.style={label={[name=#1]below:$ $}},auto]

\node{}[grow=down]
  child {node [vertex, label=30:{$\mathbf{i}$}]{}[level distance=12mm]
  child {node [vertex,label=180:{$\mathbf{k}_1$}]{}edge from parent [solid]}
  child {node [subtree, rotate=60,label=190:{$A$}]{}edge from parent [solid]}  
   child {node [vertex, label=0:{$\mathbf{j}$}]{} [level distance=12mm]
   child {node [vertex, label=0:{$\uline{\mathbf{n\!+\!2}}(u_1)$}]{}edge from parent [solid]} edge from parent[solid]} 
   child {node [subtree, rotate=-30,label=-60:{$B$}]{}edge from parent [solid]} 
   edge from parent [dashed]};
\end{tikzpicture}
        \caption{$\mathbf{n+2}$ is a singleton leaf of $\widetilde{T}$ and the parent of $\mathbf{n+2}$ has at least two elder siblings.}
        \label{(four-1)}
    \end{subfigure}\qquad
  \begin{subfigure}[t]{0.4\textwidth}
    \centering        \begin{tikzpicture}
[vertex/.style={shape=circle, draw, inner sep=1pt, fill=black},
subtree/.style={shape=ellipse, draw,minimum width=1.5cm, minimum height=.5cm},
every fit/.style={ellipse,draw,inner sep=-2pt},
sibling distance=1.5cm,level distance=1cm,
leaf/.style={label={[name=#1]below:$ $}},auto]

\node{}[grow=down]
  child {node [vertex, label=30:{$\mathbf{i}$}]{}[level distance=12mm]
  child {node [vertex,label=180:{$\mathbf{k}_1$}]{}edge from parent [solid]}
  child {node [subtree, rotate=60,label=190:{$A$}]{}edge from parent [solid]}  
   child {node [vertex, label=180:{$\uline{\mathbf{n\!+\!2}}$}]{} [level distance=12mm]
   child {node [vertex, label=0:{$\mathbf{j}(u_1)$}]{}edge from parent [solid]} edge from parent[solid]} 
   child {node [subtree, rotate=-30,label=-60:{$B$}]{}edge from parent [solid]} 
   edge from parent [dashed]};
\end{tikzpicture}
        \caption{$\mathbf{n+2}$ is the parent of singleton leaf in $\widetilde{T}$ and $\mathbf{n+2}$ has at least two elder siblings.}
        \label{(four-2)}
    \end{subfigure}
  
     \caption{$v_1\rightarrow 2u_1$.}
     \label{(four-i)}
\end{figure}
\begin{figure}[H]
    \centering
     \begin{subfigure}[t]{0.4\textwidth}
   \centering
        \begin{tikzpicture}
[vertex/.style={shape=circle, draw, inner sep=1pt, fill=black},
subtree/.style={shape=ellipse, draw,minimum width=1.5cm, minimum height=.5cm},
sibling distance=1.8cm,level distance=0.8cm,
leaf/.style={label={[name=#1]below:$#1$}},auto]
\node{}[grow=down]
  child {node [vertex, label=0:{$\mathbf{j}$}]{}[level distance=15mm]
   child {node [vertex, inner sep=1pt, label=270:{$\uline{\mathbf{n\!+\! 2}}(u_3)$}]{}edge from parent [solid] }
   child {node [vertex, label=0:{$\mathbf{i}$}]{} 
   child{node[vertex,label=270:{$\mathbf{k}_1$}] {} edge from parent [solid]}
   child {node [subtree, dashed,rotate=-55,label=-60:{$A$}]{}edge from parent [solid]} 
   edge from parent [solid]}
  child {node [subtree, rotate=-40,label=-60:{$B$}]{}edge from parent [solid]} edge from parent [dashed]
};
\end{tikzpicture}
        \caption{$\mathbf{n+2}$ is an elder non-twin leaf  of $\widetilde{T}$ and the vertex right immediately after $\mathbf{n+2}$ has at least two children.}
        \label{(four-3)}
    \end{subfigure}\qquad
   \begin{subfigure}[t]{0.4\textwidth}
   \centering
         \begin{tikzpicture}
[vertex/.style={shape=circle, draw, inner sep=1pt, fill=black},
subtree/.style={shape=ellipse, draw,minimum width=1.5cm, minimum height=.5cm},
sibling distance=1.8cm,level distance=0.8cm,
leaf/.style={label={[name=#1]below:$#1$}},auto]
\node{}[grow=down]
  child {node [vertex, label=0:{$\uline{\mathbf{n\!+\! 2}}$}]{}[level distance=15mm]
   child {node [vertex, inner sep=1pt, label=270:{$\mathbf{j}(u_3)$}]{}edge from parent [solid] }
   child {node [vertex, label=0:{$\mathbf{i}$}]{} 
   child{node[vertex,label=270:{$\mathbf{k}_1$}] {} edge from parent [solid]}
   child {node [subtree, dashed,rotate=-55,label=-60:{$A$}]{}edge from parent [solid]} 
   edge from parent [solid]}
  child {node [subtree, rotate=-40,label=-60:{$B$}]{}edge from parent [solid]} edge from parent [dashed]
};
\end{tikzpicture}
        \caption{$\mathbf{n+2}$ is  the parent of an  elder non-twin leaf in $\widetilde{T}$ and the second child of $\mathbf{n+2}$ has at least two children.}
        \label{(four-4)}
    \end{subfigure}
     \caption{$v_1\rightarrow 2u_3$.}
     \label{(four-ii)}
\end{figure}
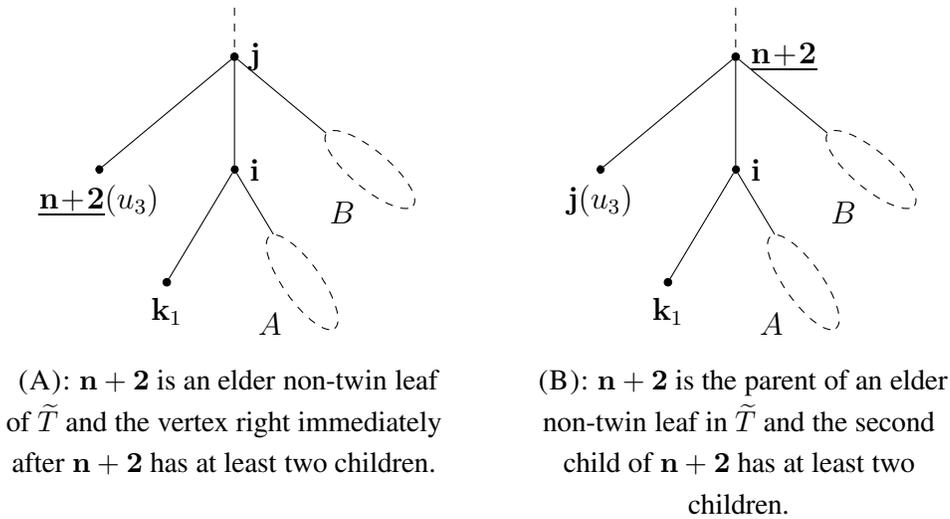

\item[(5).] The position $\diamond$  is a second leaf $\mathbf{j}$ of $T$ labeled by  $v_2$, as depicted in Fig. \ref{(five)}. 
We perform the steps shown in Fig. \ref{(five-i)} to create the  plane tree $\widetilde{T}$. It is easy to check that the change in weights follows the substitution rule $v_2\rightarrow 4u_1u_2^{-1}u_3$.

\begin{figure}[H]
    \centering
    
      \begin{tikzpicture}
[vertex/.style={shape=circle, draw, inner sep=1pt, fill=black},
subtree/.style={shape=ellipse, draw,minimum width=1.5cm, minimum height=.5cm},
sibling distance=1.8cm,level distance=0.8cm,
leaf/.style={label={[name=#1]below:$#1$}},auto]
\node{}[grow=down]
  child {node [vertex, label=0:{$\mathbf{i}$}]{}[level distance=15mm]
   child {node [vertex, inner sep=1pt, label=270:{$\mathbf{k}(u_2)$}]{}edge from parent [solid] }
   child {node [diamond,solid,draw,inner sep=1.3pt, label=270:{$\mathbf{j}(v_2)$}]{}  edge from parent [solid]}
   child {node [subtree,solid, rotate=140]{}  edge from parent [solid]}edge from parent [dashed]
};
\end{tikzpicture}
 \caption{The position $\diamond$ is a second  leaf of $T$.}
    \label{(five)}
\end{figure}
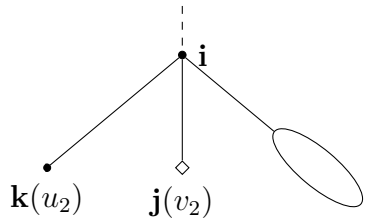

\begin{figure}[H]
    \centering
    \begin{subfigure}[t]{0.4\textwidth}
   \centering 
\begin{tikzpicture}
[vertex/.style={shape=circle, draw, inner sep=1pt, fill=black},
subtree/.style={shape=ellipse, draw,minimum width=1.5cm, minimum height=.5cm},
sibling distance=1.8cm,level distance=0.8cm,
leaf/.style={label={[name=#1]below:$#1$}},auto]
\node{}[grow=down]
  child {node [vertex, label=0:{$\mathbf{i}$}]{}[level distance=15mm]
   child {node [vertex, inner sep=1pt, label=270:{$\mathbf{k}(u_3)$}]{}edge from parent [solid] }
   child {node [vertex, label=0:{$\mathbf{j}$}]{} 
   child {node [vertex, label=0:{$\uline{ \mathbf{n\!+\!2}}(u_1)$}]{} }
   edge from parent [solid]}
   child {node [subtree,solid, rotate=140]{}  edge from parent [solid]}edge from parent [dashed]
};
\end{tikzpicture}
        \caption{$\mathbf{n+2}$ is a singleton leaf of $\widetilde{T}$ and the parent of $\mathbf{n+2}$ has only elder sibling.}
        \label{(five-1)}
    \end{subfigure}
    \begin{subfigure}[t]{0.4\textwidth}
   \centering 
\begin{tikzpicture}
[vertex/.style={shape=circle, draw, inner sep=1pt, fill=black},
subtree/.style={shape=ellipse, draw,minimum width=1.5cm, minimum height=.5cm},
sibling distance=1.8cm,level distance=0.8cm,
leaf/.style={label={[name=#1]below:$#1$}},auto]
\node{}[grow=down]
  child {node [vertex, label=0:{$\mathbf{i}$}]{}[level distance=15mm]
   child {node [vertex, inner sep=1pt, label=270:{$\mathbf{k}(u_3)$}]{}edge from parent [solid] }
   child {node [vertex, label=0:{$\uline{\mathbf{n\!+\!2}}$}]{} 
   child {node [vertex, label=0:{$ \mathbf{j}(u_1)$}]{} }
   edge from parent [solid]}
   child {node [subtree,solid, rotate=140]{}  edge from parent [solid]}edge from parent [dashed]
};
\end{tikzpicture}
        \caption{$\mathbf{n+2}$ is the parent of a  singleton leaf in $\widetilde{T}$  and  $\mathbf{n+2}$ has only elder sibling.}
        \label{(five-2)}
    \end{subfigure}
\begin{subfigure}[t]{0.4\textwidth}
\centering
     \begin{tikzpicture}
[vertex/.style={shape=circle, draw, inner sep=1pt, fill=black},
subtree/.style={shape=ellipse, draw,minimum width=1.5cm, minimum height=.5cm},
sibling distance=1.8cm,level distance=0.8cm,
leaf/.style={label={[name=#1]below:$#1$}},auto]
\node{}[grow=down]
  child {node [vertex, label=0:{$\mathbf{j}$}]{}[level distance=15mm]
   child {node [vertex, inner sep=1pt, label=270:{$\uline{\mathbf{n\!+\! 2}}(u_3)$}]{}edge from parent [solid] }
   child {node [vertex, label=0:{$\mathbf{i}$}]{} 
   child{node[vertex,label=0:{$\mathbf{k}(u_1)$}] {} edge from parent [solid]}
   edge from parent [solid]}
   child {node [subtree,solid, rotate=140]{}  edge from parent [solid]}edge from parent [dashed]
};
\end{tikzpicture}
\caption{$\mathbf{n+2}$ is an elder non-twin leaf  of $\widetilde{T}$ and the vertex right immediately after $\mathbf{n+2}$ has only one child.}
\label{(five-3)}
\end{subfigure}\qquad
\begin{subfigure}[t]{0.4\textwidth}
    \begin{tikzpicture}
[vertex/.style={shape=circle, draw, inner sep=1pt, fill=black},
subtree/.style={shape=ellipse, draw,minimum width=1.5cm, minimum height=.5cm},
sibling distance=1.8cm,level distance=0.8cm,
leaf/.style={label={[name=#1]below:$#1$}},auto]
\node{}[grow=down]
  child {node [vertex, label=0:{$ \uline{\mathbf{n\!+\! 2}}$}]{}[level distance=15mm]
   child {node [vertex, inner sep=1pt, label=270:{$\mathbf{j}(u_3)$}]{}edge from parent [solid] }
   child {node [vertex, label=0:{$\mathbf{i}$}]{} 
   child{node[vertex,label=0:{$\mathbf{k}(u_1)$}] {} edge from parent [solid]}
   edge from parent [solid]}
   child {node [subtree,solid, rotate=140]{}  edge from parent [solid]}edge from parent [dashed]
};
\end{tikzpicture}
\caption{$\mathbf{n+2}$ is the parent of an elder non-twin leaf in $\widetilde{T}$ and the second child of $\mathbf{n+2}$ has only one child.}
\label{(five-4)}
\end{subfigure}

    \caption{$v_2\rightarrow 4u_1u_2^{-1}u_3$.}
    \label{(five-i)}
\end{figure}

    \item[(6).] The position $\diamond$ is a young edge $(\mathbf{i},\mathbf{j})$ of $T$ with the labeling $t$, as shown in Fig. \ref{(six)}.  Insert the edge $(\mathbf{i},\mathbf{n+2})$ right immediately  after $(\mathbf{i},\mathbf{j})$ as illustrated in Fig. \ref{(six-i)}. The change is accordance with  $t\rightarrow t^2v_1$.

\begin{figure}[H]
    \centering
   \begin{tikzpicture}
[vertex/.style={shape=circle, draw, inner sep=1pt, fill=black},
subtree/.style={shape=ellipse, draw,minimum width=1.5cm, minimum height=.5cm},
every fit/.style={ellipse,draw,inner sep=-2pt},
sibling distance=1.8cm,level distance=1cm,
leaf/.style={label={[name=#1]below:$ $}},auto]

\node{}[grow=down]
  child {node [vertex, label=30:{$\mathbf{i}$}]{}[level distance=20mm]
  child {node [vertex]{}edge from parent [solid]}
  child {node [subtree, solid,draw, rotate=65]{}edge from parent [solid]}  
   child {node [vertex, label=0:{$\mathbf{j}$}]{} [level distance=15mm] 
   child {node [subtree, rotate= 90]{}edge from parent [solid]} edge from parent[solid] node[right,xshift=-2.5mm]{$\diamond (t)$}} 
   child {node [subtree,solid,draw, rotate=-35]{}edge from parent [solid]} 
   edge from parent [dashed]};
   
\end{tikzpicture}
    \caption{The position $\diamond$ is a young  edge.}
    \label{(six)}
\end{figure}
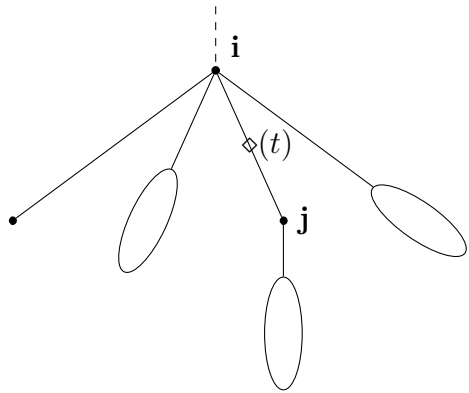

    \begin{figure}[H]
    \centering
    \begin{subfigure}[t]{0.4\textwidth}
        \centering\begin{tikzpicture}
[vertex/.style={shape=circle, draw, inner sep=1pt, fill=black},
subtree/.style={shape=ellipse, draw,minimum width=1.5cm, minimum height=.5cm},
every fit/.style={ellipse,draw,inner sep=-2pt},
sibling distance=2cm,level distance=1cm,
leaf/.style={label={[name=#1]below:$ $}},auto]

\node{}[grow=down]
  child {node [vertex, label=30:{$\mathbf{i}$}]{}[level distance=20mm]
  child {node [vertex]{}edge from parent [solid]}
  child {node [subtree, solid,draw, rotate=40]{}edge from parent [solid]}  
   child {node [vertex, label=0:{$\mathbf{j}$}]{} [level distance=15mm]
   child {node [subtree, rotate= 90]{}edge from parent [solid]} edge from parent[solid] node[right]{$(t)$}} 
   child{node[vertex,label=270:{$\uline{\mathbf{n\!+\!2}}(v_1)$}] {}edge from parent[solid] node[right]{$(t)$}}
   child {node [subtree,solid,draw, rotate=-25]{}edge from parent [solid]} 
   edge from parent [dashed]};
   
\end{tikzpicture}
\caption{$\mathbf{n+2}$ is a  younger leaf of $\widetilde{T}$.}
    \end{subfigure}
   
    \caption{$t\rightarrow t^2v_1$.}
    \label{(six-i)}
\end{figure}
\end{itemize}

From the above construction, we find that 
\begin{itemize}
     \item  $\mathbf{n+2}$ is a singleton leaf of $\widetilde{T}$ in Fig.  \ref{(four-1)} and Fig.  \ref{(five-1)}.
    \item   $\mathbf{n+2}$ is an elder twin leaf of $\widetilde{T}$  in  Fig.  \ref{(one-3)}, Fig.  \ref{(two-3)} and Fig.  \ref{(three-3)}.
     \item  $\mathbf{n+2}$ is an elder non-twin leaf of $\widetilde{T}$ in Fig.  \ref{(four-3)} and  Fig.  \ref{(five-3)}.
      \item  $\mathbf{n+2}$ is a second leaf of $\widetilde{T}$ in Fig.  \ref{(one-1)}, Fig. \ref{(two-1)} and Fig.  \ref{(three-1)}.
       \item  $\mathbf{n+2}$ is a younger leaf of $\widetilde{T}$ in Fig.  \ref{(six-i)}.
        \item  $\mathbf{n+2}$ is a parent of a singleton  leaf of $\widetilde{T}$ in Fig. \ref{(four-2)} and Fig. \ref{(five-2)}.
        \item  $\mathbf{n+2}$ is a parent of an elder twin  leaf of $\widetilde{T}$ in Fig. \ref{(one-2)}, Fig. \ref{(two-2)} and Fig. \ref{(three-2)}.
        \item  $\mathbf{n+2}$ is a parent of an elder non-twin  leaf of $\widetilde{T}$ in Fig. \ref{(four-4)} and Fig. \ref{(five-4)}.
\end{itemize}
Hence,  this construction is  reversible according to the position of the vertex $\mathbf{n+2}$ in $\widetilde{T}$. 
Therefore, we could generate 
all   tip-augmented plane trees in $\mathcal{LT}_{n+1}$ from those plane trees in $\mathcal{LT}_{n}$ based on the above construction. Moreover, the change of weights between   tip-augmented plane trees in $\mathcal{LT}_{n+1}$  and  tip-augmented plane trees in $\mathcal{LT}_{n}$  aligns with the substitution rule given in the grammar \eqref{Refi-Motz}.
Thus we show that \eqref{gram-inductaat} is valid, and so 
 \eqref{thm-Refi-Motzee} is  also valid for $n$. This completes the proof.  \qed

\subsection{The grammatical derivation for Theorem \ref{thm-gmotzkin}}   Considering the formal derivative $D$ with respect to the grammar $M$ as given in \eqref{Refi-Motz}, 
  we observe that 
$$
\begin{aligned}
    & D(t^{-1}) = -v_1,\\[5pt]
 & D^2(t^{-1}) =  -2(u_1+u_3),\\[5pt]
 & D^3(t^{-1}) =  -12tv_2u_2 = -2D(2u_3).
\end{aligned}
$$
Thus, by applying Theorem \ref{thm-Refi-Motz}, we conclude that for  $n\geq 3$,  
\[D^n(t^{-1}) = -2D^{n-2}(2u_3)=-2n!{M_{n-2}}(u_1,u_2,u_3;v_1,v_2;t).\]
Under the assumption ${M_{0}}(u_1,u_2,u_3;v_1,v_2;t)=u_3$, we derive that 
\begin{align}
 {\rm Gen}(t^{-1};q) &= t^{-1}-v_1q-(u_1+u_3)q^2+\sum_{n\geq 3}  D^n(t^{-1}) \frac{q^n}{n!} \nonumber \\[5pt]
   &= t^{-1}-v_1q-(u_1+u_3)q^2-2\sum_{n\geq 3} n!{M_{n-2}}(u_1,u_2,u_3;v_1,v_2;t)\frac{q^{n}}{n!} \nonumber \\[5pt]
    &= t^{-1}-v_1q-(u_1+u_3)q^2-2q^2\sum_{n\geq 1} {M_n}(u_1,u_2,u_3;v_1,v_2;t) q^{n} \nonumber \\[5pt]
     &= t^{-1}-v_1q-(u_1+u_3)q^2-2q^2\sum_{n\geq 0} {M_n}(u_1,u_2,u_3;v_1,v_2;t) q^{n}+2u_3q^2 \nonumber \\[5pt]
     &= t^{-1}-v_1q+(u_3-u_1)q^2-2q^2\sum_{n\geq 0} {M_n}(u_1,u_2,u_3;v_1,v_2;t) q^{n}. \label{aatttcccdd}
\end{align}
Since 
\begin{equation}\label{aatttrrcc}
{\rm Gen}(t^{-2};q)={\rm Gen}^2(t^{-1};q),
\end{equation}
it suffices to consider ${\rm Gen}(t^{-2};q)$. Observe that   
$$
\begin{aligned}
 & D(t^{-2}) = -2t^{-1}v_1,\\[5pt]
 & D^2(t^{-2}) =  2v_1^2-4t^{-1}(u_1+u_3),\\[5pt]
 & D^3(t^{-2}) = 12v_1(u_1+u_3)-24u_2v_2,  \\[5pt]
 & D^4(t^{-2}) = 24(u_1-u_3)^2, 
 \\[5pt]
 & D^5(t^{-2}) = 0. 
\end{aligned} 
$$
It follows that 
\begin{align}
   { \rm Gen}(t^{-2};q) &= \sum_{n\geq 0} D^{n}(t^{-2})\frac{q^n}{n!} \nonumber \\[5pt]
   &=t^{-2}-2t^{-1}v_1q+(v_1^2-2t^{-1}(u_1+u_3))q^2+ (2u_1v_1+2u_3v_1-4u_2v_2)q^3 \nonumber \\[5pt]
   &\quad+(u_1-u_3)^2q^4. \label{aatttccc}
\end{align}
By setting $q=0$ in \eqref{aatttrrcc}, we  infer from \eqref{aatttcccdd} and \eqref{aatttccc}  that 
\begin{equation}\label{relgft1t2}
{\rm Gen}(t^{-1};q)=\sqrt{{\rm Gen}(t^{-2};q)}.
\end{equation}
Hence we arrive at the generating function for ${M_n}(u_1,u_2,u_3;v_1,v_2;t)$ stated in  Theorem  \ref{thm-gmotzkin} by utilizing \eqref{aatttcccdd}, \eqref{aatttccc} and \eqref{relgft1t2}.  This completes the proof. \qed

 \vskip 0.2cm
\noindent{\bf Acknowledgment.} We would like to express our gratitude to Y.-D. Sun for bringing several references related to Coker's formula to our attention. This work
was supported by the National Natural Science Foundation of China.

\end{document}